\documentclass[11pt,twoside]{article}
\usepackage{amsmath, amsthm, amscd, amsfonts, amssymb, graphicx, color}

\date{}

\begin{document}

\centerline{}

\centerline {\Large{\bf Continuous controlled generalized fusion frames in Hilbert spaces }}

\newcommand{\mvec}[1]{\mbox{\bfseries\itshape #1}}
\centerline{}
\centerline{\textbf{Prasenjit Ghosh}}
\centerline{Department of Pure Mathematics, University of Calcutta,}
\centerline{35, Ballygunge Circular Road, Kolkata, 700019, West Bengal, India}
\centerline{e-mail: prasenjitpuremath@gmail.com}
\centerline{}
\centerline{\textbf{T. K. Samanta}}
\centerline{Department of Mathematics, Uluberia College,}
\centerline{Uluberia, Howrah, 711315,  West Bengal, India}
\centerline{e-mail: mumpu$_{-}$tapas5@yahoo.co.in}

\newtheorem{Theorem}{\quad Theorem}[section]

\newtheorem{definition}[Theorem]{\quad Definition}

\newtheorem{theorem}[Theorem]{\quad Theorem}

\newtheorem{remark}[Theorem]{\quad Remark}

\newtheorem{corollary}[Theorem]{\quad Corollary}

\newtheorem{note}[Theorem]{\quad Note}

\newtheorem{lemma}[Theorem]{\quad Lemma}

\newtheorem{example}[Theorem]{\quad Example}

\newtheorem{result}[Theorem]{\quad Result}
\newtheorem{conclusion}[Theorem]{\quad Conclusion}

\newtheorem{proposition}[Theorem]{\quad Proposition}

\begin{abstract}
\textbf{\emph{We introduce the notion of continuous controlled $g$-fusion frame in Hilbert space which is the generalization of discrete controlled $g$-fusion frame and give an example.\,Some characterizations of continuous controlled $g$-fusion frame have been presented.\,We define the frame operator and multiplier of continuous controlled \,$g$-fusion Bessel families in Hilbert spaces.\,Continuous resolution of the identity operator on a Hilbert space using the theory of continuous controlled $g$-fusion frame is being considered.\,Finally, we discuss perturbation results of continuous controlled $g$-fusion frame.}}
\end{abstract}
\\
\\

\section{Introduction}
 
\smallskip\hspace{.6 cm}In 1952, Duffin and Schaeffer \cite{Duffin} introduced frame for Hilbert space to study some fundamental problems in non-harmonic Fourier series.\,Later on, after some decades, frame theory was popularized by Daubechies et al.\,\cite{Daubechies}.

Frame for Hilbert space was defined as a sequence of basis-like elements in Hilbert space.\,A sequence \,$\left\{\,f_{\,i}\,\right\}_{i \,=\, 1}^{\infty} \,\subseteq\, H$\, is called a frame for a separable Hilbert space \,$\left(\,H,\,\left<\,\cdot,\,\cdot\,\right>\,\right)$, if there exist positive constants \,$0 \,<\, A \,\leq\, B \,<\, \infty$\, such that
\[ A\; \|\,f\,\|^{\,2} \,\leq\, \sum\limits_{i \,=\, 1}^{\infty}\, \left|\ \left <\,f \,,\, f_{\,i} \, \right >\,\right|^{\,2} \,\leq\, B \,\|\,f\,\|^{\,2}\; \;\text{for all}\; \;f \,\in\, H.\]
For the past few years many other types of frames were proposed such as \,$K$-frame \cite{L}, fusion frame \cite{Kutyniok}, \,$g$-frame \cite{Sun}, \,$g$-fusion frame \cite{G, Ahmadi} and \,$K$-$g$-fusion frame \cite{Sadri} etc.\,P. Ghosh and T. K. Samanta \cite{Ghosh} have discussed generalized atomic subspaces for operators in Hilbert spaces.

Controlled frame is one of the newest generalization of frame.\,P. Balaz et al.\,\cite{B} introduced controlled frame to improve the numerical efficiency of interactive algorithms for inverting the frame operator.\,In recent times, several generalizations of controlled frame namely, controlled\,$K$-frame \cite{N}, controlled\,$g$-frame \cite{F}, controlled fusion frame \cite{AK}, controlled $g$-fusion frame \cite{HS}, controlled $K$-$g$-fusion frame \cite{GR} etc. have been appeared.\,Continuous frames were proposed by Kaiser \cite{Ka} and it was independently studied by Ali et al.\,\cite{Al}.\,At present, frame theory has been widely used in signal and image processing, filter bank theory, coding and communications, system modeling and so on.

In this paper,\, continuous controlled \,$g$-fusion frame in Hilbert space is presented and some of their properties are going to be established.\,We will see that any continuous controlled $g$-fusion frame is a continuous $g$-fusion frame and converse part is also true under some sufficient conditions.\,We consider the frame operator for a pair of continuous controlled \,$g$-fusion Bessel families.\,Multiplier of continuous controlled \,$g$-fusion Bessel families in Hilbert spaces is also discussed.\,Some useful results about continuous resolution of the identity operator on a Hilbert space using the theory of continuous controlled \,$g$-fusion frame is constructed.\,At the end, we study some perturbation results of continuous controlled $g$-fusion frame.  

Throughout this paper,\;$H$\; is considered to be a separable Hilbert space with associated inner product \,$\left <\,\cdot \,,\, \cdot\,\right>$\, and \,$\mathbb{H}$\, is the collection of all closed subspace of \,$H$.\,$(\,X,\,\mu\,)$\, denotes abstract measure space with positive measure \,$\mu$.\,$I_{H}$\; is the identity operator on \,$H$.\,$\mathcal{B}\,(\,H_{\,1},\, H_{\,2}\,)$\; is a collection of all bounded linear operators from \,$H_{\,1} \,\text{to}\, H_{\,2}$.\,In particular \,$\mathcal{B}\,(\,H\,)$\, denotes the space of all bounded linear operators on \,$H$.\;For \,$S \,\in\, \mathcal{B}\,(\,H\,)$, we denote \,$\mathcal{N}\,(\,S\,)$\; and \,$\mathcal{R}\,(\,S\,)$\, for null space and range of \,$S$, respectively.\,Also, \,$P_{M} \,\in\, \mathcal{B}\,(\,H\,)$\; is the orthonormal projection onto a closed subspace \,$M \,\subset\, H$.\;$\mathcal{G}\,\mathcal{B}\,(\,H\,)$\, denotes the set of all bounded linear operators which have bounded inverse.\,If \,$S,\, R \,\in\, \mathcal{G}\,\mathcal{B}\,(\,H\,)$, then \,$R^{\,\ast},\, R^{\,-\, 1}$\, and \,$S\,R$\, are also belongs to \,$\mathcal{G}\,\mathcal{B}\,(\,H\,)$.\;$\mathcal{G}\,\mathcal{B}^{\,+}\,(\,H\,)$\, is the set of all positive operators in \,$\mathcal{G}\,\mathcal{B}\,(\,H\,)$\, and \,$T,\, U$\, are invertible operators in \,$\mathcal{G}\,\mathcal{B}\,(\,H\,)$.

\section{Preliminaries}
\smallskip\hspace{.6 cm}In this section, we recall some necessary definitions and theorems.

\begin{theorem}(\,Douglas' factorization theorem\,)\,{\cite{Douglas}}\label{th1}
Let \;$S,\, V \,\in\, \mathcal{B}\,(\,H\,)$.\,Then the following conditions are equivalent:
\begin{description}
\item[$(i)$]$\mathcal{R}\,(\,S\,) \,\subseteq\, \mathcal{R}\,(\,V\,)$.
\item[$(ii)$]\;\;$S\, S^{\,\ast} \,\leq\, \lambda^{\,2}\; V\,V^{\,\ast}$\; for some \,$\lambda \,>\, 0$.
\item[$(iii)$]$S \,=\, V\,W$\, for some bounded linear operator \,$W$\, on \,$H$.
\end{description}
\end{theorem}

\begin{theorem}\cite{O}\label{th1.001}
The set \,$\mathcal{S}\,(\,H\,)$\; of all self-adjoint operators on \,$H$\; is a partially ordered set with respect to the partial order \,$\leq$\, which is defined as for \,$R,\,S \,\in\, \mathcal{S}\,(\,H\,)$ 
\[R \,\leq\, S \,\Leftrightarrow\, \left<\,R\,f,\, f\,\right> \,\leq\, \left<\,S\,f,\, f\,\right>\; \;\forall\; f \,\in\, H.\] 
\end{theorem}

\begin{definition}\cite{Kreyzig}
A self-adjoint operator \,$U \,:\, H \,\to\, H$\, is called positive if \,$\left<\,U\,x \,,\,  x\,\right> \,\geq\, 0$\, for all \,$x \,\in\, H$.\;In notation, we can write \,$U \,\geq\, 0$.\;A self-adjoint operator \,$V \,:\, H \,\to\, H$\, is called a square root of \,$U$\, if \,$V^{\,2} \,=\, U$.\;If, in addition \,$V \,\geq\, 0$, then \,$V$\, is called positive square root of \,$U$\, and is denoted by \,$V \,=\, U^{1 \,/\, 2}$. 
\end{definition}

\begin{theorem}\cite{Kreyzig}\label{th1.05}
The positive square root \,$V \,:\, H \,\to\, H$\, of an arbitrary positive self-adjoint operator \,$U \,:\, H \,\to\, H$\, exists and is unique.\;Further, the operator \,$V$\, commutes with every bounded linear operator on \,$H$\, which commutes with \,$U$.
\end{theorem}

In a complex Hilbert space, every bounded positive operator is self-adjoint and any two bounded positive operators can be commute with each other.

\begin{theorem}\cite{Gavruta}\label{th1.01}
Let \,$M \,\subset\, H$\; be a closed subspace and \,$T \,\in\, \mathcal{B}\,(\,H\,)$.\;Then \,$P_{\,M}\, T^{\,\ast} \,=\, P_{\,M}\,T^{\,\ast}\, P_{\,\overline{T\,M}}$.\;If \,$T$\; is an unitary operator (\,i\,.\,e \,$T^{\,\ast}\, T \,=\, I_{H}$\,), then \,$P_{\,\overline{T\,M}}\;T \,=\, T\,P_{\,M}$.
\end{theorem}

\begin{definition}\cite{Ahmadi}
Let \,$\left\{\,W_{j}\,\right\}_{ j \,\in\, J}$\, be a collection of closed subspaces of \,$H$\; and \,$\left\{\,v_{j}\,\right\}_{ j \,\in\, J}$\; be a collection of positive weights, \,$\left\{\,H_{j}\,\right\}_{ j \,\in\, J}$\, be a sequence of Hilbert spaces and let \,$\Lambda_{j} \,\in\, \mathcal{B}\,(\,H,\, H_{j}\,)$\; for each \,$j \,\in\, J$.\;Then \,$\Lambda \,=\, \{\,\left(\,W_{j},\, \Lambda_{j},\, v_{j}\,\right)\,\}_{j \,\in\, J}$\; is called a generalized fusion frame or a g-fusion frame for \,$H$\; respect to \,$\left\{\,H_{j}\,\right\}_{j \,\in\, J}$\; if there exist constants \,$0 \,<\, A \,\leq\, B \,<\, \infty$\, such that
\begin{equation}\label{eq1}
A \;\left \|\,f \,\right \|^{\,2} \,\leq\, \sum\limits_{\,j \,\in\, J}\,v_{j}^{\,2}\, \left\|\,\Lambda_{j}\,P_{\,W_{j}}\,(\,f\,) \,\right\|^{\,2} \,\leq\, B \; \left\|\, f \, \right\|^{\,2}\; \;\forall\; f \,\in\, H.
\end{equation}
The constants \,$A$\; and \,$B$\; are called the lower and upper bounds of g-fusion frame, respectively.\,If \,$A \,=\, B$\; then \,$\Lambda$\; is called tight g-fusion frame and if \;$A \,=\, B \,=\, 1$\, then we say \,$\Lambda$\; is a Parseval g-fusion frame.\;If  \,$\Lambda$\; satisfies only the right inequality of (\ref{eq1}) it is called a g-fusion Bessel sequence in \,$H$\, with bound \,$B$. 
\end{definition}

Define the space
\[l^{\,2}\left(\,\left\{\,H_{j}\,\right\}_{ j \,\in\, J}\,\right) \,=\, \left \{\,\{\,f_{\,j}\,\}_{j \,\in\, J} \,:\, f_{\,j} \;\in\; H_{j},\; \sum\limits_{\,j \,\in\, J}\, \left \|\,f_{\,j}\,\right \|^{\,2} \,<\, \infty \,\right\}\]
with inner product is given by 
\[\left<\,\{\,f_{\,j}\,\}_{ j \,\in\, J} \,,\, \{\,g_{\,j}\,\}_{ j \,\in\, J}\,\right> \;=\; \sum\limits_{\,j \,\in\, J}\, \left<\,f_{\,j} \,,\, g_{\,j}\,\right>_{H_{j}}.\]\,Clearly \,$l^{\,2}\left(\,\left\{\,H_{j}\,\right\}_{ j \,\in\, J}\,\right)$\; is a Hilbert space with the pointwise operations \cite{Sadri}. 

\begin{definition}\cite{HS}
Let \,$\left\{\,W_{j}\,\right\}_{ j \,\in\, J}$\, be a collection of closed subspaces of \,$H$\, and \,$\left\{\,v_{j}\,\right\}_{ j \,\in\, J}$\, be a collection of positive weights.\,Let \,$\left\{\,H_{j}\,\right\}_{ j \,\in\, J}$\, be a sequence of Hilbert spaces, \,$T,\, U \,\in\, \mathcal{G}\,\mathcal{B}\,(\,H\,)$\, and \,$\Lambda_{j} \,\in\, \mathcal{B}\,(\,H,\, H_{j}\,)$\, for each \,$j \,\in\, J$.\,Then the family \,$\Lambda_{T\,U} \,=\, \left\{\,\left(\,W_{j},\, \Lambda_{j},\, v_{j}\,\right)\,\right\}_{j \,\in\, J}$\, is a \,$(\,T,\,U\,)$-controlled $g$-fusion frame for \,$H$\, if there exist constants \,$0 \,<\, A \,\leq\, B \,<\, \infty$\, such that 
\begin{equation}\label{eqn1.1}
A\,\|\,f\,\|^{\,2} \,\leq\, \sum\limits_{\,j \,\in\, J}\, v^{\,2}_{j}\,\left<\,\Lambda_{j}\,P_{\,W_{j}}\,U\,f,\,  \Lambda_{j}\,P_{\,W_{j}}\,T\,f\,\right> \,\leq\, \,B\,\|\,f \,\|^{\,2}\; \;\forall\; f \,\in\, H.
\end{equation}
If \,$A \,=\, B$\, then \,$\Lambda_{T\,U}$\, is called \,$(\,T,\,U\,)$-controlled tight g-fusion frame and if \,$A \,=\, B \,=\, 1$\, then we say \,$\Lambda_{T\,U}$\, is a \,$(\,T,\,U\,)$-controlled Parseval g-fusion frame.\,If \,$\Lambda_{T\,U}$\, satisfies only the right inequality of (\ref{eqn1.1}) it is called a \,$(\,T,\,U\,)$-controlled g-fusion Bessel sequence in \,$H$.    
\end{definition}

\begin{definition}\cite{HS}
Let \,$\Lambda_{T\,U}$\, be a \,$(\,T,\,U\,)$-controlled g-fusion Bessel sequence in \,$H$\, with a bound \,$B$.\,The synthesis operator \,$T_{C} \,:\, \mathcal{K}_{\,\Lambda_{j}} \,\to\, H$\, is defined as 
\[T_{C}\,\left(\,\left\{\,v_{\,j}\,\left(\,T^{\,\ast}\,P_{\,W_{j}}\, \Lambda_{j}^{\,\ast}\,\Lambda_{j}\,P_{\,W_{j}}\,U\,\right)^{1 \,/\, 2}\,f\,\right\}_{j \,\in\, J}\,\right) \,=\, \sum\limits_{\,j \,\in\, J}\,v^{\,2}_{j}\,T^{\,\ast}\,P_{\,W_{j}}\, \Lambda_{j}^{\,\ast}\,\Lambda_{j}\,P_{\,W_{j}}\,U\,f,\]for all \,$f \,\in\, H$\, and the analysis operator \,$T^{\,\ast}_{C} \,:\, H \,\to\, \mathcal{K}_{\,\Lambda_{j}}$\,is given by 
\[T_{C}^{\,\ast}\,f \,=\,  \left\{\,v_{\,j}\,\left(\,T^{\,\ast}\,P_{\,W_{j}}\, \Lambda_{j}^{\,\ast}\,\Lambda_{j}\,P_{\,W_{j}}\,U\,\right)^{1 \,/\, 2}\,f\,\right\}_{j \,\in\, J}\; \;\forall\; f \,\in\, H,\]
where 
\[\mathcal{K}_{\,\Lambda_{j}} \,=\, \left\{\,\left\{\,v_{\,j}\,\left(\,T^{\,\ast}\,P_{\,W_{j}}\, \Lambda_{j}^{\,\ast}\,\Lambda_{j}\,P_{\,W_{j}}\,U\,\right)^{1 \,/\, 2}\,f\,\right\}_{j \,\in\, J} \,:\, f \,\in\, H\,\right\} \,\subset\, l^{\,2}\left(\,\left\{\,H_{j}\,\right\}_{ j \,\in\, J}\,\right).\]
The frame operator \,$S_{C} \,:\, H \,\to\, H$\; is defined as follows:
\[S_{C}\,f \,=\, T_{C}\,T_{C}^{\,\ast}\,f \,=\, \sum\limits_{\,j \,\in\, J}\, v_{j}^{\,2}\,T^{\,\ast}\,P_{\,W_{j}}\, \Lambda_{j}^{\,\ast}\,\Lambda_{j}\,P_{\,W_{j}}\,U\,f\; \;\forall\; f \,\in\, H\]and it is easy to verify that 
\[\left<\,S_{C}\,f,\, f\,\right> \,=\, \sum\limits_{\,j \,\in\, J}\, v^{\,2}_{j}\,\left<\,\Lambda_{j}\,P_{\,W_{j}}\,U\,f,\,  \Lambda_{j}\,P_{\,W_{j}}\,T\,f\,\right>\; \;\forall\; f \,\in\, H.\]
Furthermore, if \,$\Lambda_{T\,U}$\, is a \,$(\,T,\,U\,)$-controlled g-fusion frame with bounds \,$A$\, and \,$B$\, then \,$A\,I_{\,H} \,\leq\,S_{C} \,\leq\, B\,I_{H}$.\,Hence, \,$S_{C}$\, is bounded, invertible, self-adjoint and positive linear operator.\,It is easy to verify that \,$B^{\,-1}\,I_{H} \,\leq\, S_{C}^{\,-1} \,\leq\, A^{\,-1}\,I_{H}$.
\end{definition}

\begin{definition}\cite{MF}
Let \,$F \,:\, X \,\to\, \mathbb{H}$\; be such that for each \,$h \,\in\, H$, the mapping \,$x \,\to\, P_{\,F\,(\,x\,)}\,(\,h\,)$\; is measurable (\,i.\,e. is weakly measurable\,), \,$v \,:\, X \,\to\, \mathbb{R}^{\,+}$\, be a measurable function and \,$\left\{\,K_{x}\,\right\}_{x \,\in\, X}$\, be a collection of Hilbert spaces.\;For each \,$x \,\in\, X$, suppose that \,$\,\Lambda_{x} \,\in\, \mathcal{B}\,(\,F\,(\,x\,) \,,\, K_{x}\,)$.\;Then \,$\Lambda_{F} \,=\, \left\{\,\left(\,F\,(\,x\,),\, \Lambda_{x},\, v\,(\,x\,)\,\right)\,\right\}_{x \,\in\, X}$\; is called a generalized continuous fusion frame or a gc-fusion frame for \,$H$\, with respect to \,$(\,X,\, \mu\,)$\, and \,$v$, if there exists \,$0 \,<\, A \,\leq\, B \,<\, \infty$\; such that
\[A\, \|\,h\,\|^{\,2} \,\leq\, \int\limits_{\,X}\, v^{\,2}\,(\,x\,)\, \left\|\,\Lambda_{x}\,P_{\,F\,(\,x\,)}\,(\,h\,)\,\right\|^{\,2}\,d\mu \,\leq\, B\, \|\,h\,\|^{\,2}\;\; \;\forall\, h \,\in\, H,\]where \,$P_{\,F\,(\,x\,)}$\, is the orthogonal projection onto the subspace \,$F\,(\,x\,)$.\;$\Lambda_{F}$\, is called a tight gc-fusion frame for \,$H$\, if \,$A \,=\, B$\, and Parseval if \,$A \,=\, B \,=\, 1$.\;If we have only the upper bound, we call \,$\Lambda_{F}$\, is a Bessel gc-fusion mapping for \,$H$.
\end{definition}

Let \,$K \,=\, \oplus_{x \,\in\, X}\,K_{x}$\, and \,$L^{\,2}\left(\,X,\, K\,\right)$\, be a collection of all measurable functions \,$\varphi \,:\, X \,\to\, K$\, such that for each \,$x \,\in\, X,\, \varphi\,(\,x\,) \,\in\, K_{x}$\; and \,$\int\limits_{\,X}\,\left\|\,\varphi\,(\,x\,)\,\right\|^{\,2}\,d\mu \,<\, \infty$. It can be verified that \,$L^{\,2}\left(\,X,\, K\,\right)$\; is a Hilbert space with inner product given by
\[\left<\,\phi,\, \varphi\,\right> \,=\, \int\limits_{\,X}\, \left<\,\phi\,(\,x\,),\, \varphi\,(\,x\,)\,\right>\,d\mu\] for \,$\phi,\, \varphi \,\in\, L^{\,2}\left(\,X,\, K\,\right)$.    

\begin{definition}\cite{MF}
Let \,$\Lambda_{F} \,=\, \left\{\,\left(\,F\,(\,x\,),\, \Lambda_{x},\, v\,(\,x\,)\,\right)\,\right\}_{x \,\in\, X}$\, be a Bessel gc-fusion mapping for \,$H$.\;Then the gc-fusion pre-frame operator or synthesis operator \,$T_{g\,F} \,:\, L^{\,2}\left(\,X,\, K\,\right) \,\to\, H$\, is defined by
\[\left<\,T_{g\,F}\,(\,\varphi\,),\, h\,\right> \,=\, \int\limits_{\,X}\, v\,(\,x\,)\,\left<\,P_{\,F\,(\,x\,)}\,\Lambda^{\,\ast}_{x}\,\left(\,\varphi\,(\,x\,)\,\right),\, h\,\right>\,d\mu\]where \,$\varphi \,\in\, L^{\,2}\left(\,X,\, K\,\right)$\, and \,$h \,\in\, H$.\,$T_{g\,F}$\, is a bounded linear mapping and its adjoint operator is given by  
\[T^{\,\ast}_{g\,F} \,:\, H \,\to\, L^{\,2}\left(\,X,\, K\,\right),\; T^{\,\ast}_{g\,F}\,(\,h\,) \,=\, \left\{\,v\,(\,x\,)\, \Lambda_{x}\,P_{\,F\,(\,x\,)}\,(\,h\,)\,\right\}_{x \,\in\, X},\; h \,\in\, H\]and \,$S_{g\,F} \,=\, T_{g\,F}\,T^{\,\ast}_{g\,F}$\, is called gc-fusion frame operator.
\end{definition}

For each \,$f,\, h \,\in\, H$,
\[\left<\,S_{g\,F}\,(\,f\,),\, h\,\right> \,=\, \int\limits_{\,X}\, v^{\,2}\,(\,x\,)\,\left<\,P_{F\,(\,x\,)}\,\Lambda^{\,\ast}_{x}\,\Lambda_{x}\,P_{F\,(\,x\,)}\,f,\; h\,\right>\,d\mu.\]The operator \,$S_{g\,F}$\, is bounded, self-adjoint, positive and invertible operator on \,$H$.

\begin{definition}\cite{MK} 
A sequence \,$\left\{\,T_{x} \,:\, H \,\to\, H \,:\, x \,\in\, X\,\right\}$\; is said to be a continuous resolution of the identity operator on \,$H$\; if for each \,$f,\, g \,\in\, H$, the following are hold:
\begin{itemize}
\item[(1)] $x \,\to\, \left<\,T_{x}\,f,\, g\,\right>$\; is measurable functional on \,$X$.
\item[(2)] $\left<\,f,\, g\,\right> \,=\, \int\limits_{\,X}\, \left<\,T_{x}\,f,\, g\,\right>\, d\mu(\,x\,)$
\end{itemize}
\end{definition}

\section{Continuous controlled $g$-fusion frame }

\smallskip\hspace{.6 cm}In this section, we give the continuous version of controlled $g$-fusion frame for \,$H$.\,Some of the recent results of controlled \,$g$-fusion frame are extended to continuous controlled $g$-fusion frame.

\begin{definition}
Let \,$F \,:\, X \,\to\, \mathbb{H}$\, be a mapping, \,$v \,:\, X \,\to\, \mathbb{R}^{\,+}$\, be a measurable function and \,$\left\{\,K_{x}\,\right\}_{x \,\in\, X}$\, be a collection of Hilbert spaces.\;For each \,$x \,\in\, X$, suppose that \,$\,\Lambda_{x} \,\in\, \mathcal{B}\,(\,F\,(\,x\,),\, K_{x}\,)$\, and \,$T,\, U \,\in\, \mathcal{G}\,\mathcal{B}^{\,+}\,(\,H\,)$.\,Then \,$\Lambda_{T\,U} \,=\, \left\{\,\left(\,F\,(\,x\,),\, \Lambda_{x},\, v\,(\,x\,)\,\right)\,\right\}_{x \,\in\, X}$\, is called a continuous \,$(\,T,\,U\,)$-controlled generalized fusion frame or continuous \,$(\,T,\,U\,)$-controlled $g$-fusion frame for \,$H$\, with respect to \,$(\,X,\, \mu\,)$\, and \,$v$, if
\begin{description}
\item[$(i)$]for each \,$f \,\in\, H$, the mapping \,$x \,\to\, P_{F\,(\,x\,)}\,(\,f\,)$\; is measurable (\,i.\,e. is weakly measurable\,).
\item[$(ii)$]there exist constants \,$0 \,<\, A \,\leq\, B \,<\, \infty$\, such that
\begin{equation}\label{eq1}
A\,\|\,f\,\|^{\,2} \,\leq\, \int\limits_{\,X}\,v^{\,2}\,(\,x\,)\,\left<\,\Lambda_{x}\,P_{\,F\,(\,x\,)}\,U\,f,\, \Lambda_{x}\,P_{\,F\,(\,x\,)}\,T\,f\,\right>\,d\mu_{x} \,\leq\, B\,\|\,f\,\|^{\,2},
\end{equation}
for all \,$f \,\in\, H$, where \,$P_{\,F\,(\,x\,)}$\, is the orthogonal projection onto the subspace \,$F\,(\,x\,)$.\,The constants \,$A,\,B$\, are called the frame bounds.
\end{description}
\end{definition}
Now, we consider the following cases:
\begin{description}
\item[$(I)$]If only the right inequality of (\ref{eq1}) holds then \,$\Lambda_{T\,U}$\, is called a continuous \,$(\,T,\,U\,)$-controlled $g$-fusion Bessel family for \,$H$.
\item[$(II)$]If \,$U \,=\, I_{H}$\, then \,$\Lambda_{T\,U}$\, is called a continuous \,$(\,T,\,I_{H}\,)$-controlled $g$-fusion frame for \,$H$. 
\item[$(III)$]If \,$T \,=\, U \,=\, I_{H}$\, then \,$\Lambda_{T\,U}$\, is called a continuous $g$-fusion frame for \,$H$. 
\end{description}
 
\begin{remark}
If the measure space \,$X \,=\, \mathbb{N}$\, and \,$\mu$\, is the counting measure then a continuous \,$(\,T,\,U\,)$-controlled $g$-fusion frame will be the discrete \,$(\,T,\,U\,)$-controlled $g$-fusion frame.   
\end{remark} 

\subsubsection{Example}
Let \,$H \,=\, \mathbb{R}^{\,3}$\, and \,$\left\{\,e_{\,1},\,e_{\,2},\, e_{\,3}\,\right\}$\, be an standard orthonormal basis for \,$H$.\,Consider
\[\mathcal{B} \,=\, \left\{\,x \,\in\, \mathbb{R}^{\,3} \,:\, \|\,x\,\| \,\leq\, 1\,\right\}.\]
Then it is a measure space equipped with the Lebesgue measure \,$\mu$.\,Suppose \,$\left\{\,B_{\,1},\,B_{\,2},\, B_{\,3}\,\right\}$\, is a partition of \,$\mathcal{B}$\, where \,$\mu\,(\,B_{1}\,) \,\geq\, \mu\,(\,B_{2}\,) \,\geq\, \mu\,(\,B_{3}\,) \,>\, 1$.\,Let \,$\mathbb{H} \,=\, \left\{\,W_{\,1},\,W_{\,2},\, W_{\,3}\,\right\}$, where \,$W_{1} \,=\, \overline{span}\,\left\{\,e_{\,2},\, e_{\,3}\,\right\}$, \,$W_{2} \,=\, \overline{span}\,\left\{\,e_{\,1},\, e_{\,3}\,\right\}$\, and \,$W_{3} \,=\, \overline{span}\,\left\{\,e_{\,1},\, e_{\,2}\,\right\}$.\,Define 
\[F \,:\, \mathcal{B} \,\to\, \mathbb{H}\hspace{.5cm}\text{by} \hspace{.3cm} F\,(\,x\,) \,=\, \begin{cases}
W_{1} & \text{if\;\;}\; x \,\in\, B_{1} \\ W_{2} & \text{if\;\;}\; x \,\in\, B_{2}\\ W_{3} & \text{if\;\;}\; x \,\in\, B_{3} \end{cases}\] 
and
\[v \,:\, \mathcal{B} \,\to\, [\,0\, \infty)\hspace{.5cm}\text{by} \hspace{.3cm} v\,(\,x\,) \,=\, \begin{cases}
1 & \text{if\;\;}\; x \,\in\, B_{1} \\ 2 & \text{if\;\;}\; x \,\in\, B_{2}\\ \,-\, 1 & \text{if\;\;}\; x \,\in\, B_{3}\,. \end{cases}\]
It is easy to verify that \,$F$\, and \,$v$\, are measurable functions.\,For each \,$x \,\in\, \mathcal{B}$, define the operator
\[\Lambda_{x}\,(\,f\,) \,=\, \dfrac{1}{\sqrt{\mu\,(\,B_{k}\,)}}\left<\,f,\, e_{k}\,\right>\,e_{k},\; \;f \,\in\, H,\]where \,$k$\, is such that \,$x \,\in\, \mathcal{B}_{k}$.\,Let \,$T\,\left(\,f_{\,1},\, f_{\,2},\, f_{\,3}\,\right) \,=\, \left(\,2\,f_{\,1},\, 3\,f_{\,2},\, 5\,f_{\,3}\,\right)$\, and \\$U\,\left(\,f_{\,1},\, f_{\,2},\, f_{\,3}\,\right) \,=\, \left(\,\dfrac{f_{\,1}}{2},\, \dfrac{f_{\,2}}{3},\, \dfrac{f_{\,3}}{4}\,\right)$\, be two operators on \,$H$.\,Then it is easy to verify that \,$T,\, U \,\in\, \mathcal{G}\,\mathcal{B}^{\,+}\,(\,H\,)$\, and $\,T\,U \,=\, U\,T$.\,Now, for any \,$f \,=\, \left(\,f_{\,1},\, f_{\,2},\, f_{\,3}\,\right) \,\in\, H$, we have 
\begin{align*}
&\int\limits_{\,\mathcal{B}}\,v^{\,2}\,(\,x\,)\,\left<\,\Lambda_{x}\,P_{\,F\,(\,x\,)}\,U\,f,\, \Lambda_{x}\,P_{\,F\,(\,x\,)}\,T\,f\,\right>\,d\mu_{x}\\
&=\,\sum\limits_{i \,=\, 1}^{\,3}\int\limits_{\,\mathcal{B}_{i}}\,v^{\,2}\,(\,x\,)\,\left<\,\Lambda_{x}\,P_{\,F\,(\,x\,)}\,U\,f,\, \Lambda_{x}\,P_{\,F\,(\,x\,)}\,T\,f\,\right>\,d\mu_{x}\\
&=\,f_{\,1}^{\,2} \,+\, 4\,f_{\,2}^{\,2} \,+\, \dfrac{5\,f_{\,3}^{\,2}}{4}.\\
&\Rightarrow\,\|\,f\,\|^{\,2} \,\leq\, \int\limits_{\,\mathcal{B}}\,v^{\,2}\,(\,x\,)\,\left<\,\Lambda_{x}\,P_{\,F\,(\,x\,)}\,U\,f,\, \Lambda_{x}\,P_{\,F\,(\,x\,)}\,T\,f\,\right>\,d\mu_{x} \,\leq\, 4\,\|\,f\,\|^{\,2}. 
\end{align*}
Thus, \,$\Lambda_{T\,U}$\, be a continuous \,$(\,T,\,U\,)$-controlled $g$-fusion frame for \,$\mathbb{R}^{\,3}$\, with bounds \,$1$\, and \,$4$.
 
\begin{proposition}
Let \,$\Lambda_{T\,U}$\, be a continuous \,$(\,T,\,U\,)$-controlled $g$-fusion Bessel family for \,$H$\, with bound \,$B$.\,Then there exists a unique bounded linear operator \,$S_{C} \,:\, H \,\to\,\ H$\, such that
\[\left<\,S_{C}\,f,\, g\,\right> \,=\, \int\limits_{\,X}\,v^{\,2}\,(\,x\,)\,\left<\,T^{\,\ast}\,P_{F\,(\,x\,)}\,\Lambda_{x}^{\,\ast}\,\Lambda_{x}\,P_{F\,(\,x\,)}\,U\,f,\, g\,\right>\,d\mu_{x}\; \;\forall\; f,\, g \,\in\, H.\]
Furthermore, if \,$\Lambda_{T\,U}$\, is a continuous \,$(\,T,\,U\,)$-controlled $g$-fusion frame for \,$H$\, then \,$A\,I_{H} \,\leq\, S_{C} \,\leq\, B\,I_{H}$. 
\end{proposition}

\begin{proof}
Define the mapping \,$\Psi \,:\, H \,\times\, H \,\to\, \mathbb{C}$\, by
\[\Psi\,(\,f,\, g\,) \,=\, \int\limits_{\,X}\,v^{\,2}\,(\,x\,)\,\left<\,T^{\,\ast}\,P_{F\,(\,x\,)}\,\Lambda_{x}^{\,\ast}\,\Lambda_{x}\,P_{F\,(\,x\,)}\,U\,f,\, g\,\right>\,d\mu_{x}\; \;\forall\; f,\, g \,\in\, H.\]
Then \,$\Psi$\, is a sesquilinear functional.\,Now, by Cauchy-Schwarz inequality, we have
\begin{align*}
&\left|\,\Psi\,(\,f,\, g\,)\,\right| \\
&\,=\, \left|\,\int\limits_{\,X}\,v^{\,2}\,(\,x\,)\,\left<\,\left(\,T^{\,\ast}\,P_{F\,(\,x\,)}\,\Lambda_{x}^{\,\ast}\,\Lambda_{x}\,P_{F\,(\,x\,)}\,U\,\right)^{1 \,/\, 2}\,f,\, \left(\,T^{\,\ast}\,P_{F\,(\,x\,)}\,\Lambda_{x}^{\,\ast}\,\Lambda_{x}\,P_{F\,(\,x\,)}\,U\,\right)^{1 \,/\, 2}\,g\,\right>\,d\mu_{x}\,\right| \\
&\leq\, \left(\,\int\limits_{\,X}\,v^{\,2}\,(\,x\,)\,\left\|\,\left(\,T^{\,\ast}\,P_{F\,(\,x\,)}\,\Lambda_{x}^{\,\ast}\,\Lambda_{x}\,P_{F\,(\,x\,)}\,U\,\right)^{1 \,/\, 2}\,f\,\right\|^{\,2}\,d\mu_{x}\,\right)^{1 \,/\, 2}\,\times\\
&\hspace{1.5cm}\left(\,\int\limits_{\,X}\,v^{\,2}\,(\,x\,)\,\left\|\,\left(\,T^{\,\ast}\,P_{F\,(\,x\,)}\,\Lambda_{x}^{\,\ast}\,\Lambda_{x}\,P_{F\,(\,x\,)}\,U\,\right)^{1 \,/\, 2}\,g\,\right\|^{\,2}\,d\mu_{x}\,\right)^{1 \,/\, 2}\\
&=\, \left(\,\int\limits_{\,X}\,v^{\,2}\,(\,x\,)\,\left<\,\Lambda_{x}\,P_{F\,(\,x\,)}\,U\,f,\, \Lambda_{x}\,P_{F\,(\,x\,)}\,T\,f\,\right>\,d\mu_{x}\,\right)^{1 \,/\, 2}\,\times\\
&\hspace{1.5cm}\left(\,\int\limits_{\,X}\,v^{\,2}\,(\,x\,)\,\left<\,\Lambda_{x}\,P_{F\,(\,x\,)}\,U\,g,\, \Lambda_{x}\,P_{F\,(\,x\,)}\,T\,g\,\right>\,d\mu_{x}\,\right)^{1 \,/\, 2}\\
&\leq\, B\,\|\,f\,\|\,\|\,g\,\|. 
\end{align*}
Thus, \,$\Psi$\, is a bounded sesquilinear functional with \,$\|\,\Psi\,\| \,\leq\, B$.\,Therefore, by Theorem 2.3.6 in \cite{GJ}, there exists a unique operator \,$S_{C} \,:\, H \,\to\,\ H$\, such that \,$\Psi\,(\,f,\, g\,) \,=\, \left<\,S_{C}\,f,\, g\,\right>$\, and \,$\|\,\Psi\,\| \,=\, \left\|\,S_{C}\,\right\|$.\,Thus, for each \,$f,\,g \,\in\, H$, we have
\[\left<\,S_{C}\,f,\, g\,\right> \,=\, \int\limits_{\,X}\,v^{\,2}\,(\,x\,)\,\left<\,T^{\,\ast}\,P_{F\,(\,x\,)}\,\Lambda_{x}^{\,\ast}\,\Lambda_{x}\,P_{F\,(\,x\,)}\,U\,f,\, g\,\right>\,d\mu_{x}.\]
Now, for each \,$f \,\in\, H$, we have 
\begin{align*}
\left<\,S_{C}\,f,\, f\,\right> &\,=\, \int\limits_{\,X}\,v^{\,2}\,(\,x\,)\,\left<\,\Lambda_{x}\,P_{F\,(\,x\,)}\,U\,f,\, \Lambda_{x}\,P_{F\,(\,x\,)}\,T\,f\,\right>\,d\mu_{x}\\
&=\,\int\limits_{\,X}\,v^{\,2}\,(\,x\,)\,\left\|\,\left(\,T^{\,\ast}\,P_{F\,(\,x\,)}\,\Lambda_{x}^{\,\ast}\,\Lambda_{x}\,P_{F\,(\,x\,)}\,U\,\right)^{1 \,/\, 2}\,f\,\right\|^{\,2}\,d\mu_{x}. 
\end{align*} 
This verifies that \,$S_{C}$\, is a positive operator.\,Also, it is easy to verify that \,$S_{C}$\, is a self-adjoint.\,Furthermore, if \,$\Lambda_{T\,U}$\, is a continuous \,$(\,T,\,U\,)$-controlled $g$-fusion frame for \,$H$\, then by (\ref{eq1}) it  is easy to verify that \,$A\,I_{H} \,\leq\, S_{C} \,\leq\, B\,I_{H}$.  
\end{proof}
 
\begin{theorem}
Let \,$\Lambda_{T\,U}$\, be a continuous \,$(\,T,\,U\,)$-controlled $g$-fusion Bessel family for \,$H$\, with bound \,$B$.\,Then the mapping \,$T_{C} \,:\, L^{\,2}\left(\,X,\, K\,\right) \,\to\, H$\, defined by 
\[\left<\,T_{C}\,\Phi,\, g\,\right> \,=\, \int\limits_{\,X}\,v^{\,2}\,(\,x\,)\,\left<\,T^{\,\ast}\,P_{F\,(\,x\,)}\,\Lambda_{x}^{\,\ast}\,\Lambda_{x}\,P_{F\,(\,x\,)}\,U\,f,\, g\,\right>\,d\mu_{x},\]
where for all \,$f \,\in\, H$, \,$\Phi \,=\, \left\{\,v\,(\,x\,)\,\left(\,T^{\,\ast}\,P_{F\,(\,x\,)}\,\Lambda_{x}^{\,\ast}\,\Lambda_{x}\,P_{F\,(\,x\,)}\,U\,\right)^{1 \,/\, 2}\,f\,\right\}_{x \,\in\, X}$\, and \,$g \,\in\, H$, is a linear and bounded operator with \,$\left\|\,T_{C}\,\right\| \,\leq\, \sqrt{B}$.\,Furthermore, for each \,$g \,\in\, H$, we have 
\[T_{C}^{\,\ast}\,g \,=\, \left\{\,v\,(\,x\,)\,\left(\,T^{\,\ast}\,P_{F\,(\,x\,)}\,\Lambda_{x}^{\,\ast}\,\Lambda_{x}\,P_{F\,(\,x\,)}\,U\,\right)^{1 \,/\, 2}\,g\,\right\}_{x \,\in\, X}.\] 
\end{theorem} 

\begin{proof}
For \,$\Phi \,=\, \left\{\,v\,(\,x\,)\,\left(\,T^{\,\ast}\,P_{F\,(\,x\,)}\,\Lambda_{x}^{\,\ast}\,\Lambda_{x}\,P_{F\,(\,x\,)}\,U\,\right)^{1 \,/\, 2}\,f\,\right\}_{x \,\in\, X} \,\in\, L^{\,2}\left(\,X,\, K\,\right)$,
\begin{align*}
\left\|\,T_{C}\,\Phi\,\right\| &\,=\, \sup\limits_{\|\,g\,\| \,=\, 1}\,\left|\,\left<\,T_{C}\,\Phi,\, g\,\right>\,\right|\\
& \,=\, \sup\limits_{\|\,g\,\| \,=\, 1}\,\left|\,\int\limits_{\,X}\,v^{\,2}\,(\,x\,)\,\left<\,T^{\,\ast}\,P_{F\,(\,x\,)}\,\Lambda_{x}^{\,\ast}\,\Lambda_{x}\,P_{F\,(\,x\,)}\,U\,f,\, g\,\right>\,d\mu_{x}\,\right|\\
&\leq\,\sup\limits_{\|\,g\,\| \,=\, 1}\,\left(\,\int\limits_{\,X}\,v^{\,2}\,(\,x\,)\,\left\|\,\left(\,T^{\,\ast}\,P_{F\,(\,x\,)}\,\Lambda_{x}^{\,\ast}\,\Lambda_{x}\,P_{F\,(\,x\,)}\,U\,\right)^{1 \,/\, 2}\,g\,\right\|^{\,2}\,d\mu_{x}\,\right)^{1 \,/\, 2}\,\|\,\Phi\,\|_{2}\\
&=\,\sup\limits_{\|\,g\,\| \,=\, 1}\,\left(\,\int\limits_{\,X}\,v^{\,2}\,(\,x\,)\,\left<\,\Lambda_{x}\,P_{F\,(\,x\,)}\,U\,g,\, \Lambda_{x}\,P_{F\,(\,x\,)}\,T\,g\,\right>\,d\mu_{x}\,\right)^{1 \,/\, 2}\,\|\,\Phi\,\|_{2}\\
&\leq\, \sqrt{B}\,\,\|\,\Phi\,\|_{2}.   
\end{align*}
This shows that \,$T_{C}$\, is a bounded linear operator with \,$\left\|\,T_{C}\,\right\| \,\leq\, \sqrt{B}$.\,Now, for each \,$g \,\in\, H$\, and \,$\Phi \,\in\, L^{\,2}\left(\,X,\, K\,\right)$, we have
\begin{align*}
&\left<\,\Phi,\, T^{\,\ast}_{C}\,g\,\right> \,=\, \left<\,T_{C}\,\Phi,\, g\,\right>\\
&=\,\int\limits_{\,X}\,v^{\,2}\,(\,x\,)\left<\,\left(\,T^{\,\ast}\,P_{F\,(\,x\,)}\,\Lambda_{x}^{\,\ast}\,\Lambda_{x}\,P_{F\,(\,x\,)}\,U\,\right)^{1 \,/\, 2}\,f,\,\left(\,T^{\,\ast}\,P_{F\,(\,x\,)}\,\Lambda_{x}^{\,\ast}\,\Lambda_{x}\,P_{F\,(\,x\,)}\,U\,\right)^{1 \,/\, 2}\,g\,\right>\,d\mu_{x}\\
&=\,\left<\,\Phi,\, \left\{\,v\,(\,x\,)\,\left(\,T^{\,\ast}\,P_{F\,(\,x\,)}\,\Lambda_{x}^{\,\ast}\,\Lambda_{x}\,P_{F\,(\,x\,)}\,U\,\right)^{1 \,/\, 2}\,g\,\right\}_{x \,\in\, X}\,\right>.   
\end{align*}
This completes the proof.
\end{proof} 

The operators \,$T_{C}$\, and \,$T_{C}^{\,\ast}$\, are called the synthesis operator and analysis operator of \,$\Lambda_{T\,U}$, respectively.\\

In the following proposition, we will see that it is enough to check the continuous controlled \,$g$-fusion frame condition on a dense subset \,$M$\, of \,$H$.

\begin{proposition}
Suppose that \,$(\,X,\, \mu\,)$\, is a measure space with \,$\mu$\, is \,$\sigma$-finite and \,$\Lambda_{T\,U}$\, is a continuous \,$(\,T,\,U\,)$-controlled $g$-fusion frame for a dense subset \,$M$\, of \,$H$\, having bounds \,$A$\, and \,$B$.\,Then \,$\Lambda_{T\,U}$\, is a continuous \,$(\,T,\,U\,)$-controlled $g$-fusion frame for \,$H$\, with same bounds.  
\end{proposition} 

\begin{proof}
Let \,$\left\{\,X_{n}\,\right\}_{n \,=\, 1}^{\,\infty}$\, be a sequence of disjoint measurable subsets of \,$X$\, such that \,$X \,=\, \bigcup_{n \,=\, 1}^{\,\infty}\,X_{n}$\, with \,$\mu\left(\,X_{n}\,\right) \,<\, \infty$\, for each \,$n \,\in\, \mathbb{N}$.\,Let
\[\Omega_{m} \,=\, \left\{\,x \,\in\, X \,:\, m \,\leq\, \|\,\phi\,(\,x\,)\,\| \,<\, m \,+\, 1,\, \,\forall\, \phi \,\in\, L^{\,2}\left(\,X,\, K\,\right)\,\right\},\, \,m \,\geq\, 0.\]
It is easy to verify that for each \,$m \,\geq\, 0$, \,$\Omega_{m}$\, is a measurable set and \,$X \,=\, \bigcup_{m \,=\, 0,\, n \,=\, 1}^{\,\infty}\,\left(\,X_{n} \,\cap\, \Omega_{m}\,\right)$.\,If possible suppose that \,$\Lambda_{T\,U}$\, is not a continuous \,$(\,T,\,U\,)$-controlled $g$-fusion Bessel mapping for \,$H$.\,Then there exists \,$f \,\in\, H$\, such that
\begin{align*}
&\int\limits_{\,X}\,v^{\,2}\,(\,x\,)\,\left<\,\Lambda_{x}\,P_{\,F\,(\,x\,)}\,U\,f,\, \Lambda_{x}\,P_{\,F\,(\,x\,)}\,T\,f\,\right>\,d\mu_{x} \,>\, B\,\|\,f\,\|^{\,2}\\
&\Rightarrow\,\sum\limits_{m \,=\, 0}^{\,\infty}\,\sum\limits_{n \,=\, 1}^{\,\infty}\,\int\limits_{X_{n} \,\cap\, \Omega_{m}}\,v^{\,2}\,(\,x\,)\,\left<\,\Lambda_{x}\,P_{\,F\,(\,x\,)}\,U\,f,\, \Lambda_{x}\,P_{\,F\,(\,x\,)}\,T\,f\,\right>\,d\mu_{x} \,>\, B\,\|\,f\,\|^{\,2}.
\end{align*}
Therefore, there exist finite subsets \,$I,\,J$\, such that
\begin{equation}\label{enq1}
\sum\limits_{m \,\in\, I}\,\sum\limits_{n \,\in\, J}\,\int\limits_{X_{n} \,\cap\, \Omega_{m}}\,v^{\,2}\,(\,x\,)\,\left<\,\Lambda_{x}\,P_{\,F\,(\,x\,)}\,U\,f,\, \Lambda_{x}\,P_{\,F\,(\,x\,)}\,T\,f\,\right>\,d\mu_{x} \,>\, B\,\|\,f\,\|^{\,2}.
\end{equation}
Let \,$\left\{\,f_{\,k}\,\right\}$\, be a sequence in \,$M$\, such that \,$f_{\,k} \,\to\, f$\, as \,$k \,\to\, \infty$.\,Then, we have
\[\sum\limits_{m \,\in\, I}\,\sum\limits_{n \,\in\, J}\,\int\limits_{X_{n} \,\cap\, \Omega_{m}}\,v^{\,2}\,(\,x\,)\,\left<\,\Lambda_{x}\,P_{\,F\,(\,x\,)}\,U\,f_{\,k},\, \Lambda_{x}\,P_{\,F\,(\,x\,)}\,T\,f_{\,k}\,\right>\,d\mu_{x} \,\leq\, B\,\|\,f_{\,k}\,\|^{\,2},\]
and therefore by Lebesgue's Dominated Convergence Theorem, it is a contradiction of (\ref{enq1}).\,Hence, \,$\Lambda_{T\,U}$\, is a continuous \,$(\,T,\,U\,)$-controlled $g$-fusion Bessel mapping for \,$H$.\,So, the analysis operator \,$T_{C}^{\,\ast}$\, is well-defined for \,$H$.\,Let \,$f \,\in\, H$\, be arbitrary and \,$\left\{\,f_{\,k}\,\right\}$\, be a sequence in \,$M$\, such that \,$f_{\,k} \,\to\, f$\, as \,$k \,\to\, \infty$.\,Then
\[A\,\left\|\,f_{\,k}\,\right\|^{\,2} \,\leq\, \left\|\,T_{C}^{\,\ast}\,f_{\,k}\,\right\|^{\,2}.\]
Taking \,$k \,\to\,\infty$, we get
\begin{align*}
A\,\left\|\,f\,\right\|^{\,2}& \,\leq\, \left\|\,T_{C}^{\,\ast}\,f\,\right\|^{\,2}\\
&=\,\int\limits_{\,X}\,v^{\,2}\,(\,x\,)\,\left<\,\Lambda_{x}\,P_{\,F\,(\,x\,)}\,U\,f,\, \Lambda_{x}\,P_{\,F\,(\,x\,)}\,T\,f\,\right>\,d\mu_{x}.
\end{align*} 
This completes the proof.
\end{proof}
 
Next we will see that continuous controlled \,$g$-fusion Bessel families for \,$H$\, becomes continuous controlled \,$g$-fusion frames for \,$H$\, under some sufficient conditions. 
Consider \,$G \,:\, X \,\to\, \mathbb{H}$\; be such that for each \,$h \,\in\, H$, the mapping \,$x \,\to\, P_{\,G\,(\,x\,)}\,(\,h\,)$\; is measurable and \,$w \,:\, X \,\to\, \mathbb{R}^{\,+}$\, be a measurable function.

\begin{theorem}
Let the families \,$\Lambda_{T\,U} \,=\, \left\{\,\left(\,F\,(\,x\,),\, \Lambda_{x},\, v\,(\,x\,)\,\right)\,\right\}_{x \,\in\, X}$\, and \,$\Gamma_{T\,U} \,=\, \left\{\,\left(\,G\,(\,x\,),\, \Gamma_{x},\, w\,(\,x\,)\,\right)\,\right\}_{x \,\in\, X}$\, be two continuous \,$(\,T,\,U\,)$-controlled $g$-fusion Bessel families for \,$H$\, with bounds \,$B$\, and \,$D$, respectively.\,Suppose that \,$T_{C}$\, and \,$T_{C^{\,\prime}}$\, be their synthesis operators such that \,$T_{C^{\,\prime}}\,T^{\,\ast}_{C} \,=\, I_{H}$.\,Then \,$\Lambda_{T\,U}$\, and \,$\Gamma_{T\,U}$\, are continuous \,$(\,T,\,U\,)$-controlled $g$-fusion frame for \,$H$.  
\end{theorem}

\begin{proof}
For each \,$f \,\in\, H$, we have
\begin{align*}
\|\,f\,\|^{\,4} &\,=\, \left<\,f,\, f\,\right>^{\,2} \,=\, \left<\,T^{\,\ast}_{C}\,f,\, T^{\,\ast}_{C^{\,\prime}}\,f\,\right>^{\,2} \,\leq\, \left\|\,T^{\,\ast}_{C}\,f\,\right\|^{\,2}\,\left\|\,T^{\,\ast}_{C^{\,\prime}}\,f\,\right\|^{\,2}\\
&=\,\int\limits_{\,X}\,v^{\,2}\,(\,x\,)\,\left<\,\Lambda_{x}\,P_{F\,(\,x\,)}\,U\,f,\, \Lambda_{x}\,P_{F\,(\,x\,)}\,T\,f\,\right>\,d\mu_{x}\;\times\\
&\hspace{1cm}\int\limits_{\,X}\,w^{\,2}\,(\,x\,)\,\left<\,\Gamma_{x}\,P_{G\,(\,x\,)}\,U\,f,\, \Gamma_{x}\,P_{G\,(\,x\,)}\,T\,f\,\right>\,d\mu_{x}\,\\
&\leq\, D\,\|\,f\,\|^{\,2}\,\int\limits_{\,X}\,v^{\,2}\,(\,x\,)\,\left<\,\Lambda_{x}\,P_{F\,(\,x\,)}\,U\,f,\, \Lambda_{x}\,P_{F\,(\,x\,)}\,T\,f\,\right>\,d\mu_{x}\\
&\Rightarrow\, \dfrac{1}{D}\,\|\,f\,\|^{\,2} \,\leq\, \int\limits_{\,X}\,v^{\,2}\,(\,x\,)\,\left<\,\Lambda_{x}\,P_{F\,(\,x\,)}\,U\,f,\, \Lambda_{x}\,P_{F\,(\,x\,)}\,T\,f\,\right>\,d\mu_{x}.  
\end{align*}
This shows that \,$\Lambda_{T\,U}$\, is a continuous \,$(\,T,\,U\,)$-controlled $g$-fusion frame for \,$H$\, with bounds \,$1 \,/\, D$\, and \,$B$.\,Similarly, it can be shown that \,$\Gamma_{T\,U}$\, is a continuous \,$(\,T,\,U\,)$-controlled $g$-fusion frame for \,$H$.  
\end{proof} 
 
In the next result, we construct continuous controlled \,$g$-fusion frame by using bounded linear operator. 

\begin{theorem}\label{ttt1}
Let \,$\Lambda_{T\,U}$\, be a continuous \,$(\,T,\,U\,)$-controlled $g$-fusion frame for \,$H$\, with bounds \,$A,\, B$\, and \,$V \,\in\, \mathcal{B}\,(\,H\,)$\, be an invertible operator on \,$H$\, such that \,$V^{\,\ast}$\, commutes with \,$T$, \,$U$.\,Then \,$\Gamma_{T\,U} \,=\, \left\{\,\left(\,V\,F\,(\,x\,),\, \Lambda_{x}\,P_{F\,(\,x\,)}\,V^{\,\ast},\, v\,(\,x\,)\,\right)\,\right\}_{x \,\in\, X}$\, is a continuous \,$(\,T,\,U\,)$-controlled $g$-fusion frame for \,$H$.  
\end{theorem} 
 
\begin{proof}
Since \,$P_{\,F\,(\,x\,)}\,V^{\,\ast} \,=\, P_{\,F\,(\,x\,)}\,V^{\,\ast}\,P_{\,V\,F\,(\,x\,)}$\, for all \,$x \,\in\, X$, the mapping \,$x \,\to\, P_{\,V\,F\,(\,x\,)}$\, is weakly measurable.\,Now, for each \,$f \,\in\, H$, using Theorem \ref{th1.01}, we have
\begin{align}
&\int\limits_{\,X}\,v^{\,2}\,(\,x\,)\,\left<\,\Lambda_{x}\,P_{\,F\,(\,x\,)}\,V^{\,\ast}\,P_{\,V\,F\,(\,x\,)}\,U\,f,\, \Lambda_{x}\,P_{F\,(\,x\,)}\,V^{\,\ast}\,P_{V\,F\,(\,x\,)}\,T\,f\,\right>\,d\mu_{x}\nonumber\\
&=\,\int\limits_{\,X}\,v^{\,2}\,(\,x\,)\,\left<\,\Lambda_{x}\,P_{F\,(\,x\,)}\,V^{\,\ast}\,U\,f,\, \Lambda_{x}\,P_{F\,(\,x\,)}\,V^{\,\ast}\,T\,f\,\right>\,d\mu_{x}\nonumber\\ 
&=\,\int\limits_{\,X}\,v^{\,2}\,(\,x\,)\,\left<\,\Lambda_{x}\,P_{F\,(\,x\,)}\,U\,V^{\,\ast}\,f,\, \Lambda_{x}\,P_{F\,(\,x\,)}\,T\,V^{\,\ast}\,f\,\right>\,d\mu_{x}\label{eq1.1}\\
&\leq\, B\,\left\|\,V^{\,\ast}\,f\,\right\|^{\,2} \,\leq\, B\,\|\,V\,\|^{\,2}\,\|\,f\,\|^{\,2}.\nonumber
\end{align}
On the other hand, from (\ref{eq1.1}), we get
\begin{align*}
&\int\limits_{\,X}\,v^{\,2}\,(\,x\,)\,\left<\,\Lambda_{x}\,P_{F\,(\,x\,)}\,V^{\,\ast}\,P_{V\,F\,(\,x\,)}\,U\,f,\, \Lambda_{x}\,P_{F\,(\,x\,)}\,V^{\,\ast}\,P_{V\,F\,(\,x\,)}\,T\,f\,\right>\,d\mu_{x}\\
&\hspace{1cm}\geq\, A\,\left\|\,V^{\,\ast}\,f\,\right\|^{\,2} \,\geq\, A\,\left\|\,V^{\,\,-\, 1}\,\right\|^{\,-\, 2}\,\|\,f\,\|^{\,2}\; \;\forall\; f \,\in\, H.
\end{align*}
Thus, \,$\Gamma_{T\,U}$\, is a continuous \,$(\,T,\,U\,)$-controlled $g$-fusion frame for \,$H$.\\

Furthermore, for each \,$f \,\in\, H$, using (\ref{eq1.1}), we have
\begin{align*}
&\int\limits_{\,X}\,v^{\,2}\,(\,x\,)\,\left<\,\Lambda_{x}\,P_{\,F\,(\,x\,)}\,V^{\,\ast}\,P_{\,V\,F\,(\,x\,)}\,U\,f,\, \Lambda_{x}\,P_{F\,(\,x\,)}\,V^{\,\ast}\,P_{V\,F\,(\,x\,)}\,T\,f\,\right>\,d\mu_{x}\\
&=\,\int\limits_{\,X}\,v^{\,2}\,(\,x\,)\,\left<\,\Lambda_{x}\,P_{F\,(\,x\,)}\,U\,V^{\,\ast}\,f,\, \Lambda_{x}\,P_{F\,(\,x\,)}\,T\,V^{\,\ast}\,f\,\right>\,d\mu_{x}\\
&\,=\,\left<\,S_{C}\,V^{\,\ast}\,f,\, V^{\,\ast}\,f\,\right> \,=\, \left<\,V\,S_{C}\,V^{\,\ast}\,f,\, f\,\right>, 
\end{align*}
where \,$S_{C}$\, is the corresponding frame operator for \,$\Lambda_{T\,U}$.
\end{proof}

In particular, if \,$V \,=\, S_{C}^{\,-\, 1}$\, then by the Theorem \ref{ttt1}, the family \,$\Lambda_{T\,U}^{\circ} \,=\, \left\{\,\left(\,S_{C}^{\,-\, 1}\,F\,(\,x\,),\, \Lambda_{x}\,P_{F\,(\,x\,)}\,S_{C}^{\,-\, 1},\, v\,(\,x\,)\,\right)\,\right\}_{x \,\in\, X}$\, is also a continuous \,$(\,T,\,U\,)$-controlled $g$-fusion frame for \,$H$.\,The family \,$\Lambda_{T\,U}^{\circ}$\, is called the canonical dual continuous controlled \,$g$-fusion frame of \,$\Lambda_{T\,U}$.\,It is easy to verify that the corresponding frame operator for \,$\Lambda_{T\,U}^{\circ}$\, is \,$S_{C}^{\,-\, 1}$.

A characterization of a continuous controlled $g$-fusion frame is given by in the next theorem. 

\begin{theorem}\label{th1.2}
The family \,$\Lambda_{T\,U}$\, is a continuous \,$(\,T,\,U\,)$-controlled $g$-fusion frame for \,$H$\, if and only if \,$\Lambda_{T\,U}$\, is a continuous \,$\left(\,T\,U,\, I_{H}\,\right)$-controlled $g$-fusion frame for \,$H$.    
\end{theorem}

\begin{proof}
For each \,$f \,\in\, H$, we have
\begin{align*}
&\int\limits_{\,X}\,v^{\,2}\,(\,x\,)\,\left<\,\Lambda_{x}\,P_{F\,(\,x\,)}\,U\,f,\, \Lambda_{x}\,P_{F\,(\,x\,)}\,T\,f\,\right>\,d\mu_{x} =\, \left<\,T\,S_{g\,F}\,U\,f,\, f\,\right> \,=\, \left<\,S_{g\,F}\,T\,U\,f,\, f\,\right> \\
&\hspace{1cm}=\,\int\limits_{\,X}\,v^{\,2}\,(\,x\,)\,\left<\,P_{F\,(\,x\,)}\,\Lambda_{x}^{\,\ast}\,\Lambda_{x}\,P_{F\,(\,x\,)}\,T\,U\,f,\, f\,\right>\,d\mu_{x}\\
&\hspace{1cm}=\, \int\limits_{\,X}\,v^{\,2}\,(\,x\,)\,\left<\,\Lambda_{x}\,P_{F\,(\,x\,)}\,T\,U\,f,\, \Lambda_{x}\,P_{F\,(\,x\,)}\,f\,\right>\,d\mu_{x},
\end{align*}
where 
\[\left<\,S_{g\,F}\,f,\, f\,\right> \,=\, \int\limits_{\,X}\,v^{\,2}\,(\,x\,)\,\left<\,P_{F\,(\,x\,)}\,\Lambda^{\,\ast}_{x}\,\Lambda_{x}\,P_{F\,(\,x\,)}\,f,\, f\,\right>\,d\mu_{x}.\]\,Hence, \,$\Lambda_{T\,U}$\, is continuous \,$(\,T,\,U\,)$-controlled $g$-fusion frame for \,$H$\, with bounds \,$A$\, and \,$B$\, is equivalent to:
\[A\,\|\,f\,\|^{\,2} \,\leq\, \int\limits_{\,X}\,v^{\,2}\,(\,x\,)\,\left<\,\Lambda_{x}\,P_{\,W_{x}}\,T\,U\,f,\, \Lambda_{x}\,P_{\,W_{x}}\,f\,\right>\,d\mu_{x} \,\leq\, B\,\|\,f\,\|^{\,2}\; \;\forall\; f \,\in\, H.\]  
Thus, \,$\Lambda_{T\,U}$\, is a continuous \,$\left(\,T\,U,\, I_{H}\,\right)$-controlled $g$-fusion frame for \,$H$\, with bounds \,$A$\, and \,$B$.  
\end{proof}

\begin{corollary}
The family \,$\Lambda_{T\,U}$\, is a continuous \,$(\,T,\,U\,)$-controlled $g$-fusion frame for \,$H$\, if and only if \,$\Lambda_{T\,U}$\, is a continuous \,$\left(\,(\,T\,U\,)^{1 \,/\, 2},\, (\,T\,U\,)^{1 \,/\, 2}\,\right)$-controlled $g$-fusion frame for \,$H$.   
\end{corollary}

The following theorem shows that any continuous controlled $g$-fusion frame is a continuous $g$-fusion frame and conversely any continuous $g$-fusion frame is a continuous controlled $g$-fusion frame under some conditions. 

\begin{theorem}\label{th1.3}
Let \,$T,\, U \,\in\, \mathcal{G}\,\mathcal{B}^{\,+}\,(\,H\,)$\, and \,$S_{g\,F}\,T \,=\, T\,S_{g\,F}$.\,Then \,$\Lambda_{T\,U}$\, is a continuous \,$(\,T,\,U\,)$-controlled $g$-fusion frame for \,$H$\, if and only if \,$\Lambda_{T\,U}$\, is a continuous $g$-fusion frame for \,$H$, where \,$S_{g\,F}$\, is the continuous \,$g$-fusion frame operator defined by
\[\left<\,S_{g\,F}\,f,\, f\,\right> \,=\, \int\limits_{\,X}\,v^{\,2}\,(\,x\,)\,\left<\,P_{F\,(\,x\,)}\,\Lambda^{\,\ast}_{x}\,\Lambda_{x}\,P_{F\,(\,x\,)}\,f,\, f\,\right>\,d\mu_{x},\; \;f \,\in\, H.\] 
\end{theorem}

\begin{proof}
First we suppose that \,$\Lambda_{T\,U}$\, is a continuous $g$-fusion frame for \,$H$\, with bounds \,$A$\, and \,$B$.\,Then for each \,$f \,\in\, H$, we have
\[A\,\|\,f\,\|^{\,2} \,\leq\, \int\limits_{\,X}\,v^{\,2}\,(\,x\,)\,\left\|\,\Lambda_{x}\,P_{F\,(\,x\,)}\,f\,\right\|^{\,2}\,d\mu_{x} \,\leq\, B\,\|\,f\,\|^{\,2}.\]
Now according to the Lemma 3.10 of \cite{NA}, we can deduced that 
\[m\,m^{\,\prime}\,A\,I_{H} \,\leq\, T\,S_{g\,F}\,U \,\leq\, M\,M^{\prime}\,B\,I_{H},\] where \,$m,\,m^{\,\prime}$\, and \,$M,\,M^{\prime}$\, are positive constants.\,Then for each \,$f \,\in\, H$, we have
\begin{align*}
&m\,m^{\,\prime}\,A\,\|\,f\,\|^{\,2} \,\leq\, \int\limits_{\,X}\,v^{\,2}\,(\,x\,)\,\left<\,T\,P_{F\,(\,x\,)}\,\Lambda^{\,\ast}_{x}\,\Lambda_{x}\,P_{F\,(\,x\,)}\,U\,f,\, f\,\right>\,d\mu_{x} \,\leq\, M\,M^{\,\prime}\,B\,\|\,f \,\|^{\,2}\\
&\Rightarrow\, m\,m^{\,\prime}\,A\,\|\,f\,\|^{\,2} \,\leq\, \int\limits_{\,X}\,v^{\,2}\,(\,x\,)\,\left<\,\Lambda_{x}\,P_{F\,(\,x\,)}\,U\,f,\, \Lambda_{x}\,P_{F\,(\,x\,)}\,T\,f\,\right>\,d\mu_{x} \,\leq\, M\,M^{\,\prime}\,B\,\|\,f \,\|^{\,2}.
\end{align*}
Hence, \,$\Lambda_{T\,U}$\, is a continuous \,$(\,T,\,U\,)$-controlled $g$-fusion frame for \,$H$.\\ 

Conversely, suppose that \,$\Lambda_{T\,U}$\, is a continuous \,$(\,T,\,U\,)$-controlled $g$-fusion frame for \,$H$\, with bounds \,$A$\, and \,$B$.\,Now, for each \,$f \,\in\, H$, we have
\begin{align*}
&A\,\|\,f\,\|^{\,2} \,=\, A\,\left\|\,(\,T\,U\,)^{1 \,/\, 2}\,(\,T\,U\,)^{\,-\, 1 \,/\, 2}\,f\,\right\|^{\,2}\\ 
&\leq\, \left\|\,(\,T\,U\,)^{1 \,/\, 2}\,\right\|^{\,2}\, \int\limits_{\,X}\,v^{\,2}\,(\,x\,)\,\left<\,\Lambda_{x}\,P_{F\,(\,x\,)}\,U\,(\,T\,U\,)^{\,-\, 1 \,/\, 2}\,f,\, \Lambda_{x}\,P_{F\,(\,x\,)}\,T\,(\,T\,U\,)^{\,-\, 1 \,/\, 2}\,f\,\right>\,d\mu_{x} \\
&=\, \left\|\,(\,T\,U\,)^{1 \,/\, 2}\,\right\|^{\,2}\,\int\limits_{\,X}\,v^{\,2}\,(\,x\,)\,\left<\,\Lambda_{x}\,P_{F\,(\,x\,)}\,U^{1 \,/\, 2}\,T^{\,-\, 1 \,/\, 2}\,f,\, \Lambda_{x}\,P_{F\,(\,x\,)}\,T^{1 \,/\, 2}\,U^{\,-\, 1 \,/\, 2}\,f\,\right>\,d\mu_{x}\\
&=\,\left\|\,(\,T\,U\,)^{1 \,/\, 2}\,\right\|^{\,2}\,\int\limits_{\,X}\,v^{\,2}\,(\,x\,)\,\left<\,U^{\,-\, 1 \,/\, 2}\,T^{1 \,/\, 2}\,P_{F\,(\,x\,)}\,\Lambda^{\,\ast}_{x}\,\Lambda_{x}\,P_{F\,(\,x\,)}\,U^{1 \,/\, 2}\,T^{\,-\, 1 \,/\, 2}\,f,\, f\,\right>\,d\mu_{x}\\ 
&= \left\|\,(\,T\,U\,)^{1 \,/\, 2}\,\right\|^{\,2} \left<\,U^{\,-\, 1 \,/\, 2}\,T^{1 \,/\, 2}\,S_{g\,F}\,U^{1 \,/\, 2}\,T^{\,-\, 1 \,/\, 2}\,f,\, f\,\right> = \left\|\,(\,T\,U\,)^{1 \,/\, 2}\,\right\|^{\,2} \left<\,S_{g\,F}\,f,\, f\,\right>\\  
&=\, \left\|\,(\,T\,U\,)^{1 \,/\, 2}\,\right\|^{\,2}\,\int\limits_{\,X}\,v^{\,2}\,(\,x\,)\,\left<\,P_{F\,(\,x\,)}\,\Lambda^{\,\ast}_{x}\,\Lambda_{x}\,P_{F\,(\,x\,)}\,f,\, f\,\right>\,d\mu_{x}\\ 
&\Rightarrow\, \dfrac{A}{\left\|\,(\,T\,U\,)^{1 \,/\, 2}\,\right\|^{\,2}}\,\|\,f\,\|^{\,2} \,\leq\, \int\limits_{\,X}\,v^{\,2}\,(\,x\,)\,\left\|\,\Lambda_{x}\,P_{F\,(\,x\,)}\,f\,\right\|^{\,2}\,d\mu_{x}.
\end{align*}
On the other hand, it is easy to verify that
\begin{align*}
&\int\limits_{\,X}\,v^{\,2}\,(\,x\,)\,\left\|\,\Lambda_{x}\,P_{F\,(\,x\,)}\,f\,\right\|^{\,2}\,d\mu_{x} \,=\, \left<\,(\,T\,U\,)^{\,-\, 1 \,/\, 2}\,(\,T\,U\,)^{1 \,/\, 2}\,S_{g\,F}\,f,\, f\,\right>\\
&= \left<\,(\,T\,U\,)^{1 \,/\, 2}\,S_{g\,F}\,f,\, (\,T\,U\,)^{\,-\, 1 \,/\, 2}\,f\,\right> = \left<\,S_{g\,F}\,(\,T\,U\,)\,(\,T\,U\,)^{\,-\, 1 \,/\, 2}\,f,\, (\,T\,U\,)^{\,-\, 1 \,/\, 2}\,f\,\right>\\
&=\,\left<\,T\,S_{g\,F}\,U\,(\,T\,U\,)^{\,-\, 1 \,/\, 2}\,f,\, (\,T\,U\,)^{\,-\, 1 \,/\, 2}\,f\,\right> =\left<\,S_{C}\,(\,T\,U\,)^{\,-\, 1 \,/\, 2}\,f,\, (\,T\,U\,)^{\,-\, 1 \,/\, 2}\,f\,\right>\\
&\,\leq\, B\,\left\|\,(\,T\,U\,)^{\,-\, 1 \,/\, 2}\,\right\|^{\,2}\,\|\,f\,\|^{\,2}.
\end{align*} 
Thus, \,$\Lambda_{T\,U}$\, is a continuous $g$-fusion frame for \,$H$.\,This completes the proof. 
\end{proof}

\section{Frame operator for a pair of continuous controlled $g$-fusion Bessel families}
\smallskip\hspace{.6 cm}In this section, the frame operator for a pair of continuous controlled $g$-fusion Bessel families in \,$H$\, is considered and some properties are going to be established.\,Also, we present multiplier of continuous controlled $g$-fusion Bessel families in \,$H$.\,We start this section by giving continuous resolution of the identity operator on \,$H$.

Let \,$\Lambda_{T\,U}$\, be a continuous \,$(\,T,\,U\,)$-controlled $g$-fusion frame for \,$H$\, with the corresponding frame operator \,$S_{C}$.\,Then for each \,$f,\, g \,\in\, H$, we have
\begin{align*}
\left<\,f,\, g\,\right> &\,=\, \left<\,S_{C}\,S^{\,-\, 1}_{C}\,f,\, g\,\right> \,=\, \left<\,S^{\,-\, 1}_{C}\,S_{C}\,f,\, g\,\right>\\
&=\, \int\limits_{\,X}\,v^{\,2}\,(\,x\,)\,\left<\,T^{\,\ast}\,P_{\,F\,(\,x\,)}\, \Lambda_{x}^{\,\ast}\, \Lambda_{x}\, P_{\,F\,(\,x\,)}\,U\,S^{\,-\, 1}_{C}\,f,\, g\,\right>\,d\mu_{x}\\
& \,=\, \int\limits_{\,X}\,v^{\,2}\,(\,x\,)\,\left<\,S^{\,-\, 1}_{C}\,T^{\,\ast}\,P_{\,F\,(\,x\,)}\, \Lambda_{x}^{\,\ast}\, \Lambda_{x}\, P_{\,F\,(\,x\,)}\,U\,f,\, g\,\right>\,d\mu_{x} 
\end{align*}
Thus, the families of bounded operators \,$\left\{\,T^{\,\ast}\,P_{\,F\,(\,x\,)}\, \Lambda_{x}^{\,\ast}\, \Lambda_{x}\, P_{\,F\,(\,x\,)}\,U\,S^{\,-\, 1}_{C}\,\right\}_{x \,\in\, X}$\, and \,$\left\{\,S^{\,-\, 1}_{C}T^{\,\ast}\,P_{\,F\,(\,x\,)}\, \Lambda_{x}^{\,\ast}\, \Lambda_{x}\, P_{\,F\,(\,x\,)}\,U\,\right\}_{x \,\in\, X}$\, are continuous resolution of the identity operator on \,$H$.  

\begin{theorem}
Let \,$\Lambda_{T\,U}$\, be a continuous \,$(\,T,\,U\,)$-controlled $g$-fusion frame for \,$H$\, with frame bounds \,$A,\,B$\, and \,$S_{C}$\, be its corresponding frame operator.\,Assume that \,$S^{\,-\, 1}_{C}$\, commutes with \,$T$\, and \,$U$.\,Then \,$\left\{\,v^{\,2}\,(\,x\,)\,T^{\,\ast}\,P_{\,F\,(\,x\,)}\,\Lambda_{x}^{\,\ast}\,T_{x}\,U\,\right\}_{x \,\in\, X}$\, is a continuous resolution of the identity operator on \,$H$, where \,$T_{x} \,=\, \Lambda_{x}\,P_{\,F\,(\,x\,)} \,S^{\,-\, 1}_{C},\, \,x \,\in\, X$. Furthermore, for each \,$f \,\in\, H$, we have 
\[\dfrac{A}{B^{\,2}}\,\|\,f\,\|^{\,2} \,\leq\, \int\limits_{\,X}\,v^{\,2}\,(\,x\,)\,\left<\,T_{x}\,U\,f,\, T_{x}\,T\,f\,\right>\,d\mu_{x} \,\leq\, \dfrac{B}{A^{\,2}}\,\|\,f\,\|^{\,2}.\] 
\end{theorem}

\begin{proof}
For \,$f,\, g \,\in\, H$, we have the reconstruction formula for \,$\Lambda_{T\,U}$:
\begin{align*}
\left<f,\,g \,\right>& =\, \int\limits_{\,X}\,v^{\,2}\,(\,x\,)\,\left<\,T^{\,\ast}\,P_{\,F\,(\,x\,)}\, \Lambda_{x}^{\,\ast}\, \Lambda_{x}\, P_{\,F\,(\,x\,)}\,U\,S^{\,-\, 1}_{C}\,f,\, g\,\right>\,d\mu_{x}\\
&\,=\, \int\limits_{\,X}\,v^{\,2}\,(\,x\,)\,\left<\,T^{\,\ast}\,P_{\,F\,(\,x\,)}\, \Lambda_{x}^{\,\ast}\, \Lambda_{x}\, P_{\,F\,(\,x\,)}\,S^{\,-\, 1}_{C}\,U\,f,\, g\,\right>\,d\mu_{x}\\
& \,=\, \int\limits_{\,X}\,v^{\,2}\,(\,x\,)\,\left<\,T^{\,\ast}\,P_{\,F\,(\,x\,)}\, \Lambda_{x}^{\,\ast}\,T_{x}\,U\,f,\, g\,\right>\,d\mu_{x}.
\end{align*}
Thus, \,$\left\{\,v^{\,2}\,(\,x\,)\,T^{\,\ast}\,P_{\,F\,(\,x\,)}\,\Lambda_{x}^{\,\ast}\,T_{x}\,U\,\right\}_{x \,\in\, X}$\, is a continuous resolution of the identity operator on \,$H$.\,Since \,$\Lambda_{T\,U}$\, is a continuous \,$(\,T,\,U\,)$-controlled $g$-fusion frame for \,$H$\, with frame bounds \,$A$\, and \,$B$, for each \,$f \,\in\, H$, we have
\begin{align*}
&\int\limits_{\,X}\,v^{\,2}\,(\,x\,)\,\left<\,T_{x}\,U\,f,\, T_{x}\,T\,f\,\right>\,d\mu_{x}\\
& \,=\, \int\limits_{\,X}\,v^{\,2}\,(\,x\,)\,\left<\,\Lambda_{x}\,P_{\,F\,(\,x\,)} \,S^{\,-\, 1}_{C}\,U\,f,\, \Lambda_{x}\,P_{\,F\,(\,x\,)} \,S^{\,-\, 1}_{C}\,T\,f\,\right>\,d\mu_{x}\\
&=\, \int\limits_{\,X}\,v^{\,2}\,(\,x\,)\,\left<\,\Lambda_{x}\,P_{\,F\,(\,x\,)}\,U\,S^{\,-\, 1}_{C}\,f,\, \Lambda_{x}\,P_{\,F\,(\,x\,)}\,T\,S^{\,-\, 1}_{C}\,f\,\right>\,d\mu_{x}\\
&\leq\, B\,\left\|\,S^{\,-\, 1}_{C}\,f\,\right\|^{\,2} \,\leq\, B\,\left\|\,S^{\,-\, 1}_{C}\,\right\|^{\,2}\,\|\,f\,\|^{\,2} \,\leq\, \dfrac{B}{A^{\,2}}\,\|\,f\,\|^{\,2}.  
\end{align*}
On the other hand, for each \,$f \,\in\, H$, we have 
\[\int\limits_{\,X}\,v^{\,2}\,(\,x\,)\,\left<\,T_{x}\,U\,f,\, T_{x}\,T\,f\,\right>\,d\mu_{x} \,\geq\, A\,\left\|\,S^{\,-\, 1}_{C}\,f\,\right\|^{\,2} \,\geq\, \dfrac{A}{B^{\,2}}\,\|\,f\,\|^{\,2}.\]
This completes the proof.   
\end{proof}

Next we will see that a continuous controlled $g$-fusion Bessel family becomes a continuous controlled $g$-fusion frame by using a continuous resolution of the identity operator on \,$H$.       

\begin{theorem}
Let \,$\Lambda_{T\,T}$\, be a continuous \,$(\,T,\,T\,)$-controlled $g$-fusion Bessel family in \,$H$\, with bound \,$B$.\,Then \,$\Lambda_{T\,T}$\, is a continuous \,$(\,U,\,U\,)$-controlled $g$-fusion frame for  \,$H$\, provided \,$\left\{\,v^{\,2}\,(\,x\,)\,T^{\,\ast}\,P_{\,F\,(\,x\,)}\,\Lambda_{x}^{\,\ast}\,\Lambda_{x}\,P_{\,F\,(\,x\,)}\,U\,\right\}_{x \,\in\, X}$\, is a continuous resolution of the identity operator on \,$H$.      
\end{theorem}

\begin{proof}
Since \,$\left\{\,v^{\,2}\,(\,x\,)\,T^{\,\ast}\,P_{\,F\,(\,x\,)}\,\Lambda_{x}^{\,\ast}\,\Lambda_{x}\,P_{\,F\,(\,x\,)}\,U\,\right\}_{x \,\in\, X}$\, is a continuous resolution of the identity operator on \,$H$, for \,$f,\, g \,\in\, H$, we have
\[\left<\,f,\, g\,\right> \,=\, \int\limits_{\,X}\,v^{\,2}\,(\,x\,)\,\left<\,T^{\,\ast}\,P_{\,F\,(\,x\,)}\,\Lambda_{x}^{\,\ast}\,\Lambda_{x}\,P_{\,F\,(\,x\,)}\,U\,f,\, g\,\right>\,d\mu_{x}.\]
By Cauchy-Schwartz inequality, for each \,$f \,\in\, H$, we have
\begin{align*}
\|\,f\,\|^{\,4} &\,=\, \left(\,\left<\,f,\, f\,\right>\,\right)^{\,2} \,=\, \left(\,\int\limits_{\,X}\,v^{\,2}\,(\,x\,)\,\left<\,T^{\,\ast}\,P_{\,F\,(\,x\,)}\,\Lambda_{x}^{\,\ast}\,\Lambda_{x}\,P_{\,F\,(\,x\,)}\,U\,f,\, g\,\right>\,d\mu_{x}\,\right)^{\,2}\\
&=\, \left(\,\int\limits_{\,X}\,v^{\,2}\,(\,x\,)\,\left<\,\Lambda_{x}\,P_{F\,(\,x\,)}\,U\,f,\, \Lambda_{x}\,P_{F\,(\,x\,)}\,T\,f\,\right>\,d\mu_{x}\,\right)^{\,2}\\
&\leq\, \int\limits_{\,X}\,v^{\,2}\,(\,x\,)\,\left\|\,\Lambda_{x}\,P_{F\,(\,x\,)}\,U\,f\,\right\|^{\,2}\,d\mu_{x}\,\int\limits_{\,X}\,v^{\,2}\,(\,x\,)\,\left\|\,\Lambda_{x}\,P_{F\,(\,x\,)}\,T\,f\,\right\|^{\,2}\,d\mu_{x} \\
&\leq\, B\,\|\,f\,\|^{\,2}\,\int\limits_{\,X}\,v^{\,2}\,(\,x\,)\,\left<\,\Lambda_{x}\,P_{F\,(\,x\,)}\,U\,f,\, \Lambda_{x}\,P_{F\,(\,x\,)}\,U\,f\,\right>\,d\mu_{x}\\
&\Rightarrow\,\dfrac{1}{B}\,\|\,f\,\|^{\,2} \,\leq\, \int\limits_{\,X}\,v^{\,2}\,(\,x\,)\,\left<\,\Lambda_{x}\,P_{F\,(\,x\,)}\,U\,f,\, \Lambda_{x}\,P_{F\,(\,x\,)}\,U\,f\,\right>\,d\mu_{x}.  
\end{align*}
On the other hand, for each \,$f \,\in\, H$, we have
\begin{align*}
&\int\limits_{\,X}\,v^{\,2}\,(\,x\,)\,\left<\,\Lambda_{x}\,P_{F\,(\,x\,)}\,U\,f,\, \Lambda_{x}\,P_{F\,(\,x\,)}\,U\,f\,\right>\,d\mu_{x}\\
&\hspace{.45cm} \,=\, \int\limits_{\,X}\,v^{\,2}\,(\,x\,)\,\left<\,\Lambda_{x}\,P_{F\,(\,x\,)}\,T\,T^{\,-\, 1}\,U\,f,\, \Lambda_{x}\,P_{F\,(\,x\,)}\,T\,T^{\,-\, 1}\,U\,f\,\right>\,d\mu_{x}\\
&\hspace{.45cm}\leq\, B\,\left\|\,T^{\,-\, 1}\,U\,f\,\right\|^{\,2} \,\leq\, B\,\left\|\,T^{\,-\, 1}\,\right\|^{\,2}\,\|\,U\,\|^{\,2}\,\|\,f\,\|^{\,2}.
\end{align*}
Thus, \,$\Lambda_{T\,T}$\, is a continuous \,$(\,U,\,U\,)$-controlled $g$-fusion frame for  \,$H$.\,Similarly, it can be shown that if \,$\Lambda_{T\,T}$\, is a continuous \,$(\,U,\,U\,)$-controlled $g$-fusion Bessel family in \,$H$\, then \,$\Lambda_{T\,T}$\, is also a continuous \,$(\,T,\,T\,)$-controlled $g$-fusion frame for \,$H$.   
\end{proof}

Suppose \,$G \,:\, X \,\to\, \mathbb{H}$\, be a weakly measurable function, \,$w \,:\, X \,\to\, \mathbb{R}^{\,+}$\, be a measurable function and for each \,$x \,\in\, X$, \,$\,\Gamma_{x} \,\in\, \mathcal{B}\,(\,G\,(\,x\,),\, K_{x}\,)$\, and \,$\Gamma_{U\,U}$\, denotes the family \,$\left\{\,\left(\,G\,(\,x\,),\, \Gamma_{x},\, w\,(\,x\,)\,\right)\,\right\}_{x \,\in\, X}$.\,Now, we present the frame operator for a pairs of continuous controlled $g$-fusion Bessel families. 

\begin{definition}
Let \,$\Lambda_{T\,T}$\, and \,$\Gamma_{U\,U}$\, be continuous \,$(\,T,\,T\,)$-controlled and \,$(\,U,\,U\,)$-controlled \,$g$-fusion Bessel families for \,$H$\, with bounds \,$B$\, and \,$D$, respectively.\,Then the operator \,$S_{\Lambda\,T\,\Gamma\,U} \,:\, H \,\to\, H$\, defined by
\[\left<\,S_{\Lambda\,T\,\Gamma\,U}\,f,\, g\,\right> \,=\, \int\limits_{\,X}\,v\,(\,x\,)\,w\,(\,x\,)\,\left<\,U\,P_{G\,(\,x\,)}\,\Gamma_{x}^{\,\ast}\,\Lambda_{x}\,P_{F\,(\,x\,)}\,T\,f,\, g\,\right>\,d\mu_{x},\; \;f,\, g \,\in\, H\]
is called the frame operator for the pair of continuous controlled g-fusion Bessel families \,$\Lambda_{T\,T}$\, and \,$\Gamma_{U\,U}$.
\end{definition}

\begin{theorem}
Let \,$S_{\Lambda\,T\,\Gamma\,U}$\, be the frame operator for the pair of continuous \,$(\,T,\,T\,)$-controlled and \,$(\,U,\,U\,)$-controlled g-fusion Bessel families \,$\Lambda_{T\,T}$\, and \,$\Gamma_{U\,U}$\, with bounds \,$B$\, and \,$D$, respectively.\,Then \,$S_{\Lambda\,T\,\Gamma\,U}$\, is well-defined and bounded operator with \,$\left\|\,S_{\Lambda\,T\,\Gamma\,U}\,\right\| \,\leq\, \sqrt{B\,D}$.
\end{theorem}

\begin{proof}
Let \,$f,\, g \,\in\, H$.\,Then by Cauchy-Schwartz inequality, we have
\begin{align}
&\left|\,\left<\,S_{\Lambda\,T\,\Gamma\,U}\,f,\, g\,\right>\,\right| \,=\, \left|\,\int\limits_{\,X}\,v\,(\,x\,)\,w\,(\,x\,)\,\left<\,U\,P_{G\,(\,x\,)}\,\Gamma_{x}^{\,\ast}\,\Lambda_{x}\,P_{F\,(\,x\,)}\,T\,f,\, g\,\right>\,d\mu_{x}\,\right|\nonumber\\
&\leq\, \int\limits_{\,X}\,v\,(\,x\,)\,w\,(\,x\,)\,\left|\,\left<\,\Lambda_{x}\,P_{F\,(\,x\,)}\,T\,f,\, \Gamma_{x}\,P_{G\,(\,x\,)}\,U\,g\,\right>\,\right|\,d\mu_{x}\nonumber\\
&\leq\, \int\limits_{\,X}\,v\,(\,x\,)\,w\,(\,x\,)\,\left\|\,\Lambda_{x}\,P_{F\,(\,x\,)}\,T\,f\,\right\|\,\left\|\,\Gamma_{x}\,P_{G\,(\,x\,)}\,U\,g\,\right\|\,d\mu_{x}\nonumber\\
&\leq\, \left(\,\int\limits_{\,X}\,v^{\,2}\,(\,x\,)\,\left\|\,\Lambda_{x}\,P_{F\,(\,x\,)}\,T\,f\,\right\|^{\,2}\,d\mu_{x}\,\right)^{1 \,/\, 2}\,\left(\,\int\limits_{\,X}\,w^{\,2}\,(\,x\,)\,\left\|\,\Gamma_{x}\,P_{G\,(\,x\,)}\,U\,g\,\right\|^{\,2}\,d\mu_{x}\,\right)^{1 \,/\, 2}\label{eeq1}\\
&\leq\, \sqrt{B\,D}\,\|\,f\,\|\,\|\,g\,\|.\nonumber
\end{align}
Thus, \,$S_{\Lambda\,T\,\Gamma\,U}$\, is a well-defined and bounded operator with \,$\left\|\,S_{\Lambda\,T\,\Gamma\,U}\,\right\| \,\leq\, \sqrt{B\,D}$.  
\end{proof}

In particular, for \,$T \,=\, U \,=\, I_{H}$, the operator \,$S_{\Lambda\,\Gamma} \,:\, H \,\to\, H$\, defined by
\[\left<\,S_{\Lambda\,\Gamma}\,f,\, g\,\right> \,=\, \int\limits_{\,X}\,v\,(\,x\,)\,w\,(\,x\,)\,\left<\,P_{G\,(\,x\,)}\,\Gamma_{x}^{\,\ast}\,\Lambda_{x}\,P_{F\,(\,x\,)}\,f,\, g\,\right>\,d\mu_{x},\; \;f,\, g \,\in\, H\]
is well-defined bounded operator.\,Also, for each \,$f,\, g \,\in\, H$, we have
\begin{align*}
\left<\,S_{\Lambda\,T\,\Gamma\,U}\,f,\, g\,\right> &\,=\, \int\limits_{\,X}\,v\,(\,x\,)\,w\,(\,x\,)\,\left<\,U\,P_{G\,(\,x\,)}\,\Gamma_{x}^{\,\ast}\,\Lambda_{x}\,P_{F\,(\,x\,)}\,T\,f,\, g\,\right>\,d\mu_{x}\\
&=\, \int\limits_{\,X}\,v\,(\,x\,)\,w\,(\,x\,)\,\left<\,f,\, T\,P_{F\,(\,x\,)}\,\Lambda_{x}^{\,\ast}\,\Gamma_{x}\,P_{G\,(\,x\,)}\,U\,g\,\right>\,d\mu_{x}=\, \left<\,f,\, S_{\Gamma\,U\,\Lambda\,T}\,g\,\right>
\end{align*} 
and hence \,$S^{\,\ast}_{T\,\Lambda\,\Gamma\,U} \,=\, S_{\Gamma\,U\,\Lambda\,T}$.

\begin{theorem}
Let \,$S_{\Lambda\,T\,\Gamma\,U}$\, be the frame operator for the pair of continuous \,$(\,T,\,T\,)$-controlled and \,$(\,U,\,U\,)$-controlled g-fusion Bessel families \,$\Lambda_{T\,T}$\, and \,$\Gamma_{U\,U}$\, with bounds \,$B$\, and \,$D$, respectively.\,Then the following statements are equivalent:
\begin{description}
\item[$(i)$]\,$S_{\Lambda\,T\,\Gamma\,U}$\, is bounded below.
\item[$(ii)$] There exists \,$K \,\in\, \mathcal{B}\,(\,H\,)$\; such that \,$\left\{\,T_{x}\,\right\}_{x \,\in\, X}$\, is a continuous resolution of the identity operator on \,$H$, where \,$T_{x} \,=\, v\,(\,x\,)\,w\,(\,x\,)\,K\,U\,P_{G\,(\,x\,)}\,\Gamma_{x}^{\,\ast}\,\Lambda_{x}\,P_{F\,(\,x\,)}\,T,\; x \,\in\, X$.  
\end{description} 
If one of the given conditions hold, then \,$\Lambda_{T\,T}$\, is a continuous \,$(\,T,\,T\,)$-controlled $g$-fusion frame for \,$H$.
\end{theorem}

\begin{proof}
$(\,i\,) \,\Rightarrow\, (\,ii\,)$\; Suppose that \,$S_{\Lambda\,T\,\Gamma\,U}$\, is bounded below.\;Then for each \,$f \,\in\, H$, there exists \,$A \,>\, 0$\; such that
\begin{align*}
\|\,f\,\|^{\,2} \,\leq\, A\;\left\|\,S_{\Lambda\,T\,\Gamma\,U}\,f\,\right\|^{\,2} &\,\Rightarrow\, \left<\,I_{H}\,f,\, f\,\right> \,\leq\, A\, \left<\,S_{\Lambda\,T\,\Gamma\,U}^{\,\ast}\,S_{\Lambda\,T\,\Gamma\,U}\,f,\, f\,\right>\\
& \,\Rightarrow\,I^{\,\ast}_{H}\,I_{\,H} \,\leq\, A\, S_{\Lambda\,T\,\Gamma\,U}^{\,\ast}\,S_{\Lambda\,T\,\Gamma\,U}.
\end{align*}  
So, by Theorem \ref{th1}, there exists \,$K \,\in\, \mathcal{B}\,(\,H\,)$\, such that \,$K\,S_{\Lambda\,T\,\Gamma\,U} \,=\, I_{H}$.\\Therefore, for each \,$f,\, g \,\in\, H$, we have
\begin{align*}
\left<\,f,\, g\,\right> \,=\, \left<\,K\,S_{\Lambda\,T\,\Gamma\,U}\,f,\, g\,\right> \,=\, \int\limits_{\,X}\,v\,(\,x\,)\,w\,(\,x\,)\,\left<\,K\,U\,P_{G\,(\,x\,)}\,\Gamma_{x}^{\,\ast}\,\Lambda_{x}\,P_{F\,(\,x\,)}\,T\,f,\, g\,\right>\,d\mu_{x}.
\end{align*}
Thus, \,$\left\{\,T_{x}\,\right\}_{x \,\in\, X}$\, is a continuous resolution of the identity operator on \,$H$, where \,$T_{x} \,=\, v\,(\,x\,)\,w\,(\,x\,)\,K\,U\,P_{G\,(\,x\,)}\,\Gamma_{x}^{\,\ast}\,\Lambda_{x}\,P_{F\,(\,x\,)}\,T,\; x \,\in\, X$. \\\\
$(\,ii\,) \,\Rightarrow\, (\,i\,)$\, Since \,$\left\{\,T_{x}\,\right\}_{x \,\in\, X}$\, is a continuous resolution of the identity operator on \,$H$, for each \,$f,\,g \,\in\, H$, we have
\[\left<\,f,\, g\,\right> \,=\, \int\limits_{\,X}\,v\,(\,x\,)\,w\,(\,x\,)\,\left<\,K\,U\,P_{G\,(\,x\,)}\,\Gamma_{x}^{\,\ast}\,\Lambda_{x}\,P_{F\,(\,x\,)}\,T\,f,\, g\,\right>\,d\mu_{x} \,=\, \left<\,K\,S_{\Lambda\,T\,\Gamma\,U}\,f,\, g\,\right>.\]Thus, \,$I_{H} \,=\, K\,S_{\Lambda\,T\,\Gamma\,U}$.\,So, by Theorem \ref{th1}, there exists some \,$\alpha \,>\, 0$\, such that \,$I_{H}\, I_{H}^{\,\ast} \,\leq\, \alpha\, S_{\Lambda\,T\,\Gamma\,U}\,S_{\Lambda\,T\,\Gamma\,U}^{\,\ast}$\, and hence \,$S_{\Lambda\,T\,\Gamma\,U}$\, is bounded below.\\

Last part:
First we suppose that \,$S_{\Lambda\,T\,\Gamma\,U}$\, is bounded below.\,Then for all \\$f \,\in\, H$, there exists \,$M \,>\, 0$\; such that \,$\left\|\,S_{\Lambda\,T\,\Gamma\,U}\,f\,\right\| \,\geq\, M\; \|\,f\,\|$\, and therefore by (\ref{eeq1}), we have 
\begin{align*}
&M^{\,2}\, \|\,f\,\|^{\,2} \,\leq\, \left\|\,S_{\Lambda\,T\,\Gamma\,U}\,f\,\right\|^{\,2} \leq\, D\,\left(\,\int\limits_{\,X}\,v^{\,2}\,(\,x\,)\,\left\|\,\Lambda_{x}\,P_{F\,(\,x\,)}\,T\,f\,\right\|^{\,2}\,d\mu_{x}\,\right)^{1 \,/\, 2}\\
&\Rightarrow\, \dfrac{M^{\,2}}{D}\, \|\,f\,\|^{\,2} \leq\, \int\limits_{\,X}\,v^{\,2}\,(\,x\,)\,\left<\,\Lambda_{x}\,P_{F\,(\,x\,)}\,T\,f,\, \Lambda_{x}\,P_{F\,(\,x\,)}\,T\,f\,\right>\,d\mu_{x}  
\end{align*}
Hence, \,$\Lambda_{T\,T}$\, is a continuous \,$(\,T,\,T\,)$-controlled \,$g$-fusion frame for \,$H$\, with bounds \,$M^{\,2} \,/\, D$\, and \,$B$.\,Similarly, it can be shown that \,$\Gamma_{U\,U}$\, is a continuous \,$(\,U,\,U\,)$-controlled \,$g$-fusion frame for \,$H$\, with bounds \,$M^{\,2} \,/\, B$\, and \,$D$.\\ 

Next, we suppose that the given condition $(ii)$ holds.\,Then for each \,$f,\, g \,\in\, H$, we have
\[\left<\,f,\,g\,\right> \,=\, \int\limits_{\,X}\,v\,(\,x\,)\,w\,(\,x\,)\,\left<\,K\,U\,P_{G\,(\,x\,)}\,\Gamma_{x}^{\,\ast}\,\Lambda_{x}\,P_{F\,(\,x\,)}\,T\,f,\, g\,\right>\,d\mu_{x},\; K \,\in\, \mathcal{B}\,(\,H\,).\]By Cauchy-Schwarz inequality, for each \,$f  \,\in\, H$, we have
\begin{align*}
&\|\,f\,\|^{\,2} \,=\, \left<\,f,\, f\,\right> \,=\, \int\limits_{\,X}\,v\,(\,x\,)\,w\,(\,x\,)\,\left<\,K\,U\,P_{G\,(\,x\,)}\,\Gamma_{x}^{\,\ast}\,\Lambda_{x}\,P_{F\,(\,x\,)}\,T\,f,\, f\,\right>\,d\mu_{x}\\
&=\,\int\limits_{\,X}\,v\,(\,x\,)\,w\,(\,x\,)\,\left<\,\Lambda_{x}\,P_{F\,(\,x\,)}\,T\,f,\, \Gamma_{x}\,P_{G\,(\,x\,)}\,U\,K^{\,\ast}\,f\,\right>\,d\mu_{x}\\
&\leq\left(\,\int\limits_{\,X}\,v^{\,2}\,(\,x\,)\,\left\|\,\Lambda_{x}\,P_{F\,(\,x\,)}\,T\,f\,\right\|^{\,2}\,d\mu_{x}\,\right)^{1 \,/\, 2}\,\left(\,\int\limits_{\,X}\,w^{\,2}\,(\,x\,)\,\left\|\,\Gamma_{x}\,P_{G\,(\,x\,)}\,U\,K^{\,\ast}\,f\,\right\|^{\,2}\,d\mu_{x}\,\right)^{1 \,/\, 2}\\
&\leq\,\sqrt{D}\,\left\|\,K^{\,\ast}\,f\right\|\,\left(\,\int\limits_{\,X}\,v^{\,2}\,(\,x\,)\,\left\|\,\Lambda_{x}\,P_{F\,(\,x\,)}\,T\,f\,\right\|^{\,2}\,d\mu_{x}\,\right)^{1 \,/\, 2}\\
&\Rightarrow\, \dfrac{1}{D\,\|\,K\,\|^{\,2}}\, \|\,f\,\|^{\,2} \,\leq\, \int\limits_{\,X}\,v^{\,2}\,(\,x\,)\,\left<\,\Lambda_{x}\,P_{F\,(\,x\,)}\,T\,f,\, \Lambda_{x}\,P_{F\,(\,x\,)}\,T\,f\,\right>\,d\mu_{x} 
\end{align*}
Therefore, in this case \,$\Lambda_{T\,T}$\, is also a continuous \,$(\,T,\,T\,)$-controlled \,$g$-fusion frame for \,$H$.       
\end{proof}

\begin{theorem}
Let \,$\Lambda_{T\,T}$\, and \,$\Gamma_{U\,U}$\, be continuous \,$(\,T,\,T\,)$-controlled and \,$(\,U,\,U\,)$-controlled \,$g$-fusion frames for \,$H$\, and \,$T,\, U,\, S_{\Lambda\,\Gamma} \,+\, S_{\Gamma\,\Lambda} \,\in\, \mathcal{G}\,\mathcal{B}^{\,+}\,(\,H\,)$\, such that they are commutes with each others.\,Then \,$S_{\Lambda\,T\,\Gamma\,U} \,+\, S_{\Gamma\,U\,\Lambda\,T}$\, is a positive operator.
\end{theorem} 
 
\begin{proof}
For each \,$f,\, g \,\in\, H$, we have
\begin{align*}
&\left<\,\left(\,S_{\Lambda\,T\,\Gamma\,U} \,+\, S_{\Gamma\,U\,\Lambda\,T}\,\right)\,f,\, g\,\right>\\
& = \int\limits_{\,X}v^{\,2}\,(\,x\,)\left<\,U\,P_{G\,(\,x\,)}\,\Gamma_{x}^{\,\ast}\,\Lambda_{x}\,P_{F\,(\,x\,)}\,T\,f,\, g\,\right>d\mu_{x}\,+ \\
&\hspace{1.5cm}\int\limits_{\,X}v^{\,2}\,(\,x\,)\left<\,T\,P_{F\,(\,x\,)}\,\Lambda_{x}^{\,\ast}\,\Gamma_{x}\,P_{G\,(\,x\,)}\,U\,f,\, g\,\right>d\mu_{x}\\
&=\,\left<\,U\,S_{\Lambda\,\Gamma}\,T\,f,\, g\,\right> \,+\, \left<\,T\,S_{\Gamma\,\Lambda}\,U\,f,\, g\,\right> =\,\left<\,U\,S_{\Lambda\,\Gamma}\,T\,f,\, g\,\right> \,+\, \left<\,U\,S_{\Gamma\,\Lambda}\,T\,f,\, g\,\right>\\
&=\,\left<\,U\,\left(\,S_{\Lambda\,\Gamma} \,+\, S_{\Gamma\,\Lambda}\,\right)\,T\,f,\, g\,\right>.  
\end{align*}
This shows that \,$S_{\Lambda\,T\,\Gamma\,U} \,+\, S_{\Gamma\,U\,\Lambda\,T} \,=\, U\,\left(\,S_{\Lambda\,\Gamma} \,+\, S_{\Gamma\,\Lambda}\,\right)\,T$.\,Since \,$T,\,U$\, and \,$S_{\Lambda\,\Gamma} \,+\, S_{\Gamma\,\Lambda}$\, are positive and commutes with each other.\,Therefore, \,$S_{\Lambda\,T\,\Gamma\,U} \,+\, S_{\Gamma\,U\,\Lambda\,T}$\, is a positive operator. 
\end{proof} 

\begin{theorem}\label{tthm1}
Let \,$\Lambda_{T\,T}$\, and \,$\Gamma_{U\,U}$\, be continuous \,$(\,T,\,T\,)$-controlled and \,$(\,U,\,U\,)$-controlled \,$g$-fusion Bessel families for \,$H$\, with bounds \,$B$\, and \,$D$, respectively.\,Let \,$m \,\in\, L^{\,\infty}\,(\,X,\,\mu\,)$.\,Then the operator \,$M_{m,\,\Lambda\,T,\,\Gamma\,U} \,:\, H \,\to\, H$\, defined by
\[\left<\,M_{m,\,\Lambda\,T,\,\Gamma\,U}\,f,\, g\,\right> \,=\, \int\limits_{\,X}\,m\,(\,x\,)\,v\,(\,x\,)\,w\,(\,x\,)\,\left<\,T\,P_{F\,(\,x\,)}\,\Lambda_{x}^{\,\ast}\,\Gamma_{x}\,P_{G\,(\,x\,)}\,U\,f,\, g\,\right>\,d\mu_{x},\]
for \,$f,\, g \,\in\, H$, is well-defined and bounded operator.
\end{theorem}

\begin{proof}
For each \,$f,\, g \,\in\, H$, we have
\begin{align*}
&\left|\,\left<\,M_{m,\,\Lambda\,T,\,\Gamma\,U}\,f,\, g\,\right>\,\right| \,=\, \left|\,\int\limits_{\,X}\,m\,(\,x\,)\,v\,(\,x\,)\,w\,(\,x\,)\,\left<\,T\,P_{F\,(\,x\,)}\,\Lambda_{x}^{\,\ast}\,\Gamma_{x}\,P_{G\,(\,x\,)}\,U\,f,\, g\,\right>\,d\mu_{x}\,\right|\\
&\leq\, \int\limits_{\,X}\,|\,m\,(\,x\,)\,|\,v\,(\,x\,)\,w\,(\,x\,)\,\left\|\,\Lambda_{x}\,P_{F\,(\,x\,)}\,T\,g\,\right\|\,\left\|\,\Gamma_{x}\,P_{G\,(\,x\,)}\,U\,f\,\right\|\,d\mu_{x}\\
&\leq \|\,m\,\|_{\infty}\left(\int\limits_{\,X}v^{\,2}\,(\,x\,)\left\|\,\Lambda_{x}\,P_{F\,(\,x\,)}\,T\,g\,\right\|^{\,2}d\mu_{x}\right)^{1 \,/\, 2}\left(\int\limits_{\,X}w^{\,2}\,(\,x\,)\left\|\,\Gamma_{x}\,P_{G\,(\,x\,)}\,U\,f\,\right\|^{\,2}d\mu_{x}\right)^{1 \,/\, 2}\\
&\leq\, \|\,m\,\|_{\infty}\,\sqrt{B\,D}\,\|\,f\,\|\,\|\,g\,\|.
\end{align*}
Thus, \,$M_{m,\,\Lambda\,T,\,\Gamma\,U}$\, is a well-defined and bounded operator with \,$\left\|\,M_{m,\,\Lambda\,T,\,\Gamma\,U}\,\right\| \,\leq\, \|\,m\,\|_{\infty}\,\sqrt{B\,D}$.  
\end{proof}

Now, multiplier of continuous controlled \,$g$-fusion Bessel families in Hilbert spaces is presented.

\begin{definition}
Let \,$\Lambda_{T\,T}$\, and \,$\Gamma_{U\,U}$\, be continuous \,$(\,T,\,T\,)$-controlled and \,$(\,U,\,U\,)$-controlled \,$g$-fusion Bessel families for \,$H$\, with bounds \,$B$\, and \,$D$, respectively.\,Let \,$m \,\in\, L^{\,\infty}\,(\,X,\,\mu\,)$.\,Then the operator \,$M_{m,\,\Lambda\,T,\,\Gamma\,U} \,:\, H \,\to\, H$\, defined by
\[\left<\,M_{m,\,\Lambda\,T,\,\Gamma\,U}\,f,\, g\,\right> \,=\, \int\limits_{\,X}\,m\,(\,x\,)\,v\,(\,x\,)\,w\,(\,x\,)\,\left<\,T\,P_{F\,(\,x\,)}\,\Lambda_{x}^{\,\ast}\,\Gamma_{x}\,P_{G\,(\,x\,)}\,U\,f,\, g\,\right>\,d\mu_{x},\]
for \,$f,\, g \,\in\, H$, is called the continuous \,$(\,T,\,U\,)$-controlled \,$g$-fusion Bessel multiplier of \,$\Lambda_{T\,T}$, \,$\Gamma_{U\,U}$\, and \,$m$. 
\end{definition} 
 
For each \,$f,\, g \,\in\, H$, we have
\begin{align*}
\left<\,M_{m,\,\Lambda\,T,\,\Gamma\,U}\,f,\, g\,\right>& \,=\, \int\limits_{\,X}\,m\,(\,x\,)\,v\,(\,x\,)\,w\,(\,x\,)\,\left<\,T\,P_{F\,(\,x\,)}\,\Lambda_{x}^{\,\ast}\,\Gamma_{x}\,P_{G\,(\,x\,)}\,U\,f,\, g\,\right>\,d\mu_{x}\\
&=\,\int\limits_{\,X}\,m\,(\,x\,)\,v\,(\,x\,)\,w\,(\,x\,)\,\left<\,f,\, U\,P_{G\,(\,x\,)}\,\Gamma_{x}^{\,\ast}\,\Lambda_{x}\,P_{F\,(\,x\,)}\,T\,g\,\right>\,d\mu_{x}\\
&=\, \left<\,f,\, M_{m,\,\Gamma\,U,\,\Lambda\,T}\,g\,\right>
\end{align*} 
and hence \,$M_{m,\,\Lambda\,T,\,\Gamma\,U}^{\,\ast} \,=\, M_{m,\,\Gamma\,U,\,\Lambda\,T}$.

\begin{theorem}
Let \,$M_{m,\,\Lambda\,T,\,\Gamma\,U}$\, be the continuous \,$(\,T,\,U\,)$-controlled \,$g$-fusion Bessel multiplier of \,$\Lambda_{T\,T}$, \,$\Gamma_{U\,U}$\, and \,$m$.\,Assume \,$\lambda \,\in\, (\,0,\,1\,)$\, such that
\[\left\|\,f \,-\, M_{m,\,\Lambda\,T,\,\Gamma\,U}\,f\,\right\| \,\leq\, \lambda\, \|\,f\,\|\;  \;\;\forall\; f \,\in\, H.\]
Then \,$\Lambda_{T\,T}$\,  and \,$\Gamma_{U\,U}$\, are continuous \,$(\,T,\,T\,)$-controlled and \,$(\,U,\,U\,)$-controlled $g$-fusion frame for \,$H$.   
\end{theorem}

\begin{proof}
For each \,$f \,\in\, H$, we have
\begin{align*}
&\left(\,1 \,-\, \lambda\,\right)\,\|\,f\,\| \,\leq\, \left\|\,M_{m,\,\Lambda\,T,\,\Gamma\,U}\,f \,\right\| \,=\, \sup\limits_{\|\,g\,\| \,=\, 1}\,\left<\,M_{m,\,\Lambda\,T,\,\Gamma\,U}\,f,\, g\,\right>\\
&=\, \sup\limits_{\|\,g\,\| \,=\, 1}\,\int\limits_{\,X}\,m\,(\,x\,)\,v\,(\,x\,)\,w\,(\,x\,)\,\left<\,T\,P_{F\,(\,x\,)}\,\Lambda_{x}^{\,\ast}\,\Gamma_{x}\,P_{G\,(\,x\,)}\,U\,f,\, g\,\right>\,d\mu_{x}\\
&\leq \sup\limits_{\|\,g\,\| \,=\, 1}\,\|\,m\,\|_{\infty}\left(\int\limits_{\,X}v^{\,2}\,(\,x\,)\left\|\,\Lambda_{x}\,P_{F\,(\,x\,)}\,T\,g\,\right\|^{\,2}d\mu_{x}\right)^{1 \,/\, 2}\,\times\\
&\hspace{1.5cm}\left(\int\limits_{\,X}w^{\,2}\,(\,x\,)\left\|\,\Gamma_{x}\,P_{G\,(\,x\,)}\,U\,f\,\right\|^{\,2}d\mu_{x}\right)^{1 \,/\, 2}\\
&\leq\, \|\,m\,\|_{\infty}\,\sqrt{B}\,\left(\int\limits_{\,X}\,w^{\,2}\,(\,x\,)\left\|\,\Gamma_{x}\,P_{G\,(\,x\,)}\,U\,f\,\right\|^{\,2}d\mu_{x}\right)^{1 \,/\, 2}\\
&\Rightarrow\, \dfrac{\left(\,1 \,-\, \lambda\,\right)^{\,2}}{B\,\|\,m\,\|_{\infty}^{\,2}}\,\|\,f\,\|^{\,2} \,\leq\, \int\limits_{\,X}w^{\,2}\,(\,x\,)\,\left<\,\Gamma_{x}\,P_{G\,(\,x\,)}\,U\,f,\, \Gamma_{x}\,P_{G\,(\,x\,)}\,U\,f\,\right>\,d\mu_{x}.
\end{align*} 
Thus, \,$\Gamma_{U\,U}$\, is a continuous \,$(\,U,\,U\,)$-controlled $g$-fusion frame for \,$H$.\,Similarly, it can be shown that \,$\Lambda_{T\,T}$\, is a continuous \,$(\,T,\,T\,)$-controlled $g$-fusion frame for \,$H$.\\ 
\end{proof}

\section{Perturbation of continuous controlled $g$-fusion frame}

\smallskip\hspace{.6 cm}In frame theory, one of the most important problem is the stability of frame under some perturbation.\,P. Casazza and Chirstensen \cite{CC} have been generalized the Paley-Wiener perturbation theorem to perturbation of frame in Hilbert space.\,P. Ghosh and T. K. Samanta \cite{P} discussed stability of dual \,$g$-fusion frame in Hilbert space.\,Like the perturbation of discrete frames, we present a perturbation of continuous controlled $g$-fusion frame. \\

\begin{theorem}
Let \,$\Lambda_{T\,U}$\, be a continuous \,$(\,T,\,U\,)$-controlled \,$g$-fusion frame for \,$H$\, with bounds \,$A,\,B$\, and \,$\Gamma_{T\,U} \,=\, \left\{\,\left(\,G\,(\,x\,),\, \Gamma_{x},\, v\,(\,x\,)\,\right)\,\right\}_{x \,\in\, X}$.\,If there exist constants \,$\lambda_{\,1},\, \lambda_{\,2},\, \mu$\, with 
\[0 \,\leq\, \lambda_{\,1},\, \lambda_{\,2} \,<\, 1,\, \, \,A\,\left(\,1 \,-\, \lambda_{\,1}\,\right) \,-\, \mu\,\int\limits_{\,X}\,v^{\,2}\,(\,x\,)\,d\mu_{x} \,>\, 0\] such that for each \,$f \,\in\, H$,
\begin{align*}
 0 &\,\leq\, \left<\,T^{\,\ast}\,\left(\,P_{G\,(\,x\,)}\,\Gamma^{\,\ast}_{x}\,\Gamma_{x}\,P_{G\,(\,x\,)} \,-\, P_{F\,(\,x\,)}\,\Lambda^{\,\ast}_{x}\,\Lambda_{x}\,P_{F\,(\,x\,)}\,\right)\,U\,f,\, f\,\right>\\ 
&\leq\,\lambda_{\,1}\,\left<\,T^{\,\ast}\,P_{F\,(\,x\,)}\,\Lambda^{\,\ast}_{x}\,\Lambda_{x}\,P_{F\,(\,x\,)}\,U\,f,\, f\,\right> \,+\, \lambda_{\,2}\,\left<\,T^{\,\ast}\,P_{G\,(\,x\,)}\,\Gamma^{\,\ast}_{x}\,\Gamma_{x}\,P_{G\,(\,x\,)}\,U\,f,\, f\,\right> \,+\, \mu\,\|\,f\,\|^{\,2} 
\end{align*} 
then \,$\Gamma_{T\,U}$\, is a continuous \,$(\,T,\,U\,)$-controlled \,$g$-fusion frame for \,$H$.   
\end{theorem} 

\begin{proof} 
For each \,$f \,\in\, H$, we have
\begin{align*}
&\int\limits_{\,X}\,v^{\,2}\,(\,x\,)\,\left<\,T^{\,\ast}\,P_{G\,(\,x\,)}\,\Gamma^{\,\ast}_{x}\,\Gamma_{x}\,P_{G\,(\,x\,)}\,U\,f,\, f\,\right>\,d\mu_{x}\\
& =\, \int\limits_{\,X}\,v^{\,2}\,(\,x\,)\,\left<\,T^{\,\ast}\,\left(\,P_{G\,(\,x\,)}\,\Gamma^{\,\ast}_{x}\,\Gamma_{x}\,P_{G\,(\,x\,)} \,-\, P_{F\,(\,x\,)}\,\Lambda^{\,\ast}_{x}\,\Lambda_{x}\,P_{F\,(\,x\,)}\,\right)\,U\,f,\, f\,\right>\,d\mu_{x} \,+\, \\
&\hspace{1cm}+\,\int\limits_{\,X}\,v^{\,2}\,(\,x\,)\,\left<\,T^{\,\ast}\,P_{F\,(\,x\,)}\,\Lambda^{\,\ast}_{x}\,\Lambda_{x}\,P_{F\,(\,x\,)}\,U\,f,\, f\,\right>\,d\mu_{x} \\
&\leq\,\left(\,1 \,+\, \lambda_{\,1}\,\right)\,\int\limits_{\,X}\,v^{\,2}\,(\,x\,)\,\left<\,T^{\,\ast}\,P_{F\,(\,x\,)}\,\Lambda^{\,\ast}_{x}\,\Lambda_{x}\,P_{F\,(\,x\,)}\,U\,f,\, f\,\right>\,d\mu_{x} \,+\, \mu\,\left\|\,f\,\right\|^{\,2}\,\int\limits_{\,X}\,v^{\,2}\,(\,x\,)\,d\mu_{x}\\
&\hspace{1cm}+\,\lambda_{2}\,\int\limits_{\,X}\,v^{\,2}\,(\,x\,)\,\left<\,T^{\,\ast}\,P_{G\,(\,x\,)}\,\Gamma^{\,\ast}_{x}\,\Gamma_{x}\,P_{G\,(\,x\,)}\,U\,f,\, f\,\right>\,d\mu_{x}\\ 
&\Rightarrow\, \left(\,1 \,-\, \lambda_{\,2}\,\right)\,\int\limits_{\,X}\,v^{\,2}\,(\,x\,)\,\left<\,T^{\,\ast}\,P_{G\,(\,x\,)}\,\Gamma^{\,\ast}_{x}\,\Gamma_{x}\,P_{G\,(\,x\,)}\,U\,f,\, f\,\right>\,d\mu_{x}\\ 
&\leq\, \left(\,1 \,+\, \lambda_{\,1}\,\right)\,\int\limits_{\,X}\,v^{\,2}\,(\,x\,)\,\left<\,T^{\,\ast}\,P_{F\,(\,x\,)}\,\Lambda^{\,\ast}_{x}\,\Lambda_{x}\,P_{F\,(\,x\,)}\,U\,f,\, f\,\right>\,d\mu_{x} \,+\, \mu\,\left\|\,f\,\right\|^{\,2}\,\int\limits_{\,X}\,v^{\,2}\,(\,x\,)\,d\mu_{x} .\\
&\Rightarrow\,\int\limits_{\,X}\,v^{\,2}\,(\,x\,)\,\left<\,\Gamma_{x}\,P_{G\,(\,x\,)}\,U\,f,\, \Gamma_{x}\,P_{G\,(\,x\,)}\,T\,f\,\right>\,d\mu_{x} \,\leq\, \left[\,\dfrac{\left(\,1 \,+\, \lambda_{\,1}\,\right)\,B \,+\, \mu\,\int\limits_{\,X}\,v^{\,2}\,(\,x\,)\,d\mu_{x}}{\left(\,1 \,-\, \lambda_{\,2}\,\right)}\,\right]\,\|\,f\,\|^{\,2}.
\end{align*} 
On the other hand, for each \,$f \,\in\, H$, we have
\begin{align*}
&\int\limits_{\,X}\,v^{\,2}\,(\,x\,)\,\left<\,T^{\,\ast}\,P_{G\,(\,x\,)}\,\Gamma^{\,\ast}_{x}\,\Gamma_{x}\,P_{G\,(\,x\,)}\,U\,f,\, f\,\right>\,d\mu_{x}\\
& \,\geq\, \int\limits_{\,X}\,v^{\,2}\,(\,x\,)\,\left<\,T^{\,\ast}\,P_{F\,(\,x\,)}\,\Lambda^{\,\ast}_{x}\,\Lambda_{x}\,P_{F\,(\,x\,)}\,U\,f,\, f\,\right>\,d\mu_{x}\,-\\
&\,-\, \int\limits_{\,X}\,v^{\,2}\,(\,x\,)\,\left<\,T^{\,\ast}\,\left(\,P_{G\,(\,x\,)}\,\Gamma^{\,\ast}_{x}\,\Gamma_{x}\,P_{G\,(\,x\,)} \,-\, P_{F\,(\,x\,)}\,\Lambda^{\,\ast}_{x}\,\Lambda_{x}\,P_{F\,(\,x\,)}\,\right)\,U\,f,\, f\,\right>\,d\mu_{x}.\\
&\Rightarrow\,\left(\,1 \,+\, \lambda_{\,2}\,\right)\,\int\limits_{\,X}\,v^{\,2}\,(\,x\,)\,\left<\,T^{\,\ast}\,P_{G\,(\,x\,)}\,\Gamma^{\,\ast}_{x}\,\Gamma_{x}\,P_{G\,(\,x\,)}\,U\,f,\, f\,\right>\,d\mu_{x}\\
&\geq\,\left(\,1 \,-\, \lambda_{\,1}\,\right)\,\int\limits_{\,X}\,v^{\,2}\,(\,x\,)\,\left<\,T^{\,\ast}\,P_{F\,(\,x\,)}\,\Lambda^{\,\ast}_{x}\,\Lambda_{x}\,P_{F\,(\,x\,)}\,U\,f,\, f\,\right>\,d\mu_{x} \,-\, \mu\,\left\|\,f\,\right\|^{\,2}\,\int\limits_{\,X}\,v^{\,2}\,(\,x\,)\,d\mu_{x}.\\
&\Rightarrow\,\int\limits_{\,X}\,v^{\,2}\,(\,x\,)\,\left<\,\Gamma_{x}\,P_{G\,(\,x\,)}\,U\,f,\, \Gamma_{x}\,P_{G\,(\,x\,)}\,T\,f\,\right>\,d\mu_{x} \,\geq\, \left[\,\dfrac{\left(\,1 \,-\, \lambda_{\,1}\,\right)\,A \,-\, \mu\,\int\limits_{\,X}\,v^{\,2}\,(\,x\,)\,d\mu_{x}}{\left(\,1 \,+\, \lambda_{\,2}\,\right)}\,\right]\,\left\|\,f\,\right\|^{\,2}.      
\end{align*}
Thus, \,$\Gamma_{T\,U}$\, is a continuous \,$(\,T,\,U\,)$-controlled \,$g$-fusion frame for \,$H$.  
\end{proof}

\begin{corollary}
Let \,$\Lambda_{T\,U}$\, be a continuous \,$(\,T,\,U\,)$-controlled \,$g$-fusion frame for \,$H$\, with bounds \,$A,\,B$\, and \,$\Gamma_{T\,U} \,=\, \left\{\,\left(\,G\,(\,x\,),\, \Gamma_{x},\, v\,(\,x\,)\,\right)\,\right\}_{x \,\in\, X}$.\,If there exists constant \,$0 \,<\, D\,\int\limits_{\,X}\,v^{\,2}\,(\,x\,)\,d\mu_{x} \,<\, A$\, such that for each \,$f \,\in\, H$,
\[0 \,\leq\, \left<\,T^{\,\ast}\,\left(\,P_{G\,(\,x\,)}\,\Gamma^{\,\ast}_{x}\,\Gamma_{x}\,P_{G\,(\,x\,)} \,-\, P_{F\,(\,x\,)}\,\Lambda^{\,\ast}_{x}\,\Lambda_{x}\,P_{F\,(\,x\,)}\,\right)\,U\,f,\, f\,\right> \,\leq\, D\,\left\|\,f\,\right\|^{\,2}\]
then \,$\Gamma_{T\,U}$\, is a continuous \,$(\,T,\,U\,)$-controlled \,$g$-fusion frame for \,$H$.
\end{corollary}

\begin{proof}
For each \,$f \,\in\, H$, we have
\begin{align*}
&\int\limits_{\,X}\,v^{\,2}\,(\,x\,)\,\left<\,\Gamma_{x}\,P_{G\,(\,x\,)}\,U\,f,\, \Gamma_{x}\,P_{G\,(\,x\,)}\,T\,f\,\right>\,d\mu_{x} \,=\, \int\limits_{\,X}\,v^{\,2}\,(\,x\,)\,\left<\,T^{\,\ast}\,P_{G\,(\,x\,)}\,\Gamma^{\,\ast}_{x}\,\Gamma_{x}\,P_{G\,(\,x\,)}\,U\,f,\, f\,\right>\,d\mu_{x}\\
& =\, \int\limits_{\,X}\,v^{\,2}\,(\,x\,)\,\left<\,T^{\,\ast}\,\left(\,P_{G\,(\,x\,)}\,\Gamma^{\,\ast}_{x}\,\Gamma_{x}\,P_{G\,(\,x\,)} \,-\, P_{F\,(\,x\,)}\,\Lambda^{\,\ast}_{x}\,\Lambda_{x}\,P_{F\,(\,x\,)}\,\right)\,U\,f,\, f\,\right>\,d\mu_{x} \,+\, \\
&\hspace{1cm}+\,\int\limits_{\,X}\,v^{\,2}\,(\,x\,)\,\left<\,T^{\,\ast}\,P_{F\,(\,x\,)}\,\Lambda^{\,\ast}_{x}\,\Lambda_{x}\,P_{F\,(\,x\,)}\,U\,f,\, f\,\right>\,d\mu_{x} \\
&\leq\,  \left(\,B \,+\, D\,\int\limits_{\,X}\,v^{\,2}\,(\,x\,)\,d\mu_{x}\,\right)\,\|\,f\,\|^{\,2}.
\end{align*}
On the other hand,
\begin{align*}
&\int\limits_{\,X}\,v^{\,2}\,(\,x\,)\,\left<\,T^{\,\ast}\,P_{G\,(\,x\,)}\,\Gamma^{\,\ast}_{x}\,\Gamma_{x}\,P_{G\,(\,x\,)}\,U\,f,\, f\,\right>\,d\mu_{x}\\
&\hspace{1cm} \,\geq\, \int\limits_{\,X}\,v^{\,2}\,(\,x\,)\,\left<\,T^{\,\ast}\,P_{F\,(\,x\,)}\,\Lambda^{\,\ast}_{x}\,\Lambda_{x}\,P_{F\,(\,x\,)}\,U\,f,\, f\,\right>\,d\mu_{x}\,-\\
&\hspace{1cm}\,-\, \int\limits_{\,X}\,v^{\,2}\,(\,x\,)\,\left<\,T^{\,\ast}\,\left(\,P_{G\,(\,x\,)}\,\Gamma^{\,\ast}_{x}\,\Gamma_{x}\,P_{G\,(\,x\,)} \,-\, P_{F\,(\,x\,)}\,\Lambda^{\,\ast}_{x}\,\Lambda_{x}\,P_{F\,(\,x\,)}\,\right)\,U\,f,\, f\,\right>\,d\mu_{x}\\
&\hspace{1cm}\geq\, \left(\,A \,-\, D\,\int\limits_{\,X}\,v^{\,2}\,(\,x\,)\,d\mu_{x}\,\right)\,\left\|\,f\,\right\|^{\,2}\; \;\forall\; f \,\in\, H.
\end{align*}
This completes the proof.
\end{proof}


\begin{thebibliography}{0} 

\bibitem{Sadri}R. Ahmadi, G. Rahimlou, V. Sadri and R. Zarghami Farfar,
\emph{Constructions of K-g fusion frames and their duals in Hilbert spaces}, Bull. Transilvania Un. Brasov, 13(62), No. 1, (2020), 17-32.

\bibitem{Al}S. T. Ali, J. P. Antonie and J. P. Gazeau,
\emph{continuous frames in Hilbert spaces}, Annals of Physics 222, (1993), 1-37.

\bibitem{NA}N. Assila, S. Kabbaj and B. Moalige,
\emph{Controlled\,$K$-fusion frame for Hilbert space}, arXiv: 2007.05110v1.

\bibitem{B}P. Balazs, J. P. Antonie and A. Grybos, \emph{Weighted and controlled frames: mutual relationship and first numerical properties,} Int. J. Wavelets, Multiresolution Info. Proc., 14 (2010), No. 1, 109-132.

\bibitem{Kutyniok}
P. Casazza and G. Kutyniok,
\emph{Frames of subspaces}, Cotemporary Math, AMS 345 (2004), 87-114.

\bibitem{CC}P. Casazza and O. Christensen,
\emph{Perturbation of operators and applications to frame theory}, J. Fourier Anal. Appl., 3 (1997), 543-557.

\bibitem{O}O. Christensen,
\emph{An introduction to frames and Riesz bases}, Birkhauser (2008).

\bibitem{Daubechies}
I. Daubechies, A. Grossmann and Y. Mayer,
\emph{Painless nonorthogonal expansions}, Journal of Mathematical Physics 27 (5) (1986) 1271-1283.

\bibitem{Douglas}R. G. Douglas,
\emph{On majorization, factorization, and range inclusion of operators on Hilbert space}. Proc. Am. Math. Soc. 17, 413-415 (1966).

\bibitem{Duffin}R. J. Duffin and A. C. Schaeffer,
\emph{A class of nonharmonic Fourier series}, Trans. Amer. Math. Soc., 72, (1952), 341-366.

\bibitem{MF}M. H. Faroughi, A. Rahimi and R. Ahmadi,
\emph{GC-fusion frames}, Methods of Functional Analysis and Topology, Vol. 16 (2010), no. 2, pp. 112-119.

\bibitem{Gavruta}P. Gavruta, 
\emph{On the duality of fusion frames}, J. Math. Anal. Appl. 333 (2007) 871-879.

\bibitem{L}
L. Gavruta,
\emph{Frames for operator}, Appl. Comput. Harmon. Anal. 32 (1), 139-144 (2012).

\bibitem{P}P. Ghosh and T. K. Samanta,
\emph{Stability of dual g-fusion frame in Hilbert spaces}, Methods of Functional Analysis and Topology, Vol. 26, no. 3, pp. 227-240.

\bibitem{Ghosh}P. Ghosh and T. K. Samanta,
\emph{Generalized atomic subspaces for operators in Hilbert spaces}, Mathematica Bohemica, Accepted.

\bibitem{G}P. Ghosh and T. K. Samanta,
\emph{Generalized fusion frame in tensor product of Hilbert spaces}, Journal of the Indian Mathematical Society, Accepted.

\bibitem{Ka}G. Kaiser,
\emph{A Friendly Guide to Wavelets}, Birkhauser (1994). 

\bibitem{MK}M. Khayyami and A. Nazari,
\emph{Construction of continuous g-frames and continuous fusion frames}, Sahand Communications in Mathematical Analysis (SCMA) Vol. 4 No. 1 (2016), 43-55.

\bibitem{AK}A. Khosravi and K. Musazadeh,
\emph{Controlled fusion frames}, Methods Funct. Anal. Topol. 18 (3), 256-265.

\bibitem{Kreyzig}E. Kreyzig,
\emph{Introductory Functional Analysis with Applications.}
Wiley, New York (1989).

\bibitem{GJ}G. J. Murphy, \emph{$C^{\,\ast}$ Algebras and Operator Theory}, Academic Press, San Diego, 1990.

\bibitem{N}M. Nouri, A. Rahimi and Sh. Najafizadeh,
\emph{Controlled $K$-frames in Hilbert spaces,} Int. J. Anal. Appl. 4 (2015), No. 2, 39-50.

\bibitem{F}A. Rahimi and A. Fereydooni,
\emph{Controlled \,$g$-frames and their \,$g$-multipliers in Hilbert spaces}, Analele Stiintifice Ale Universitatii Ovidius Constanta Seria Matematica 21 (2), 223-236.

\bibitem{GR}G. Rahimlou, V. Sadri and R. Ahmadi,
\emph{Construction of controlled $K$-$g$-fusion frame in Hilbert spaces}, U. P. B. Sci. Bull., Series A, Vol. 82, Iss. 1, 2020.

\bibitem{Ahmadi}V. Sadri, Gh. Rahimlou, R. Ahmadi and R. Zarghami Farfar,
\emph{Generalized Fusion Frames in Hilbert Spaces}, Submitted, arXiv: 1806.03598v1  [\,math.FA\,] 10 Jun 2018. 

\bibitem{HS}H. Shakoory, R. Ahamadi, N. Behzadi and S. Nami,
\emph{$(\,C,\, C^{\,\prime}\,)$-Controlled $g$-fusion frames,} Submitted (2018).

\bibitem{Sun}W. Sun,
\emph{G-frames and G-Riesz bases}, Journal of Mathematical Analysis and Applications 322 (1) (2006), 437-452.


\end{thebibliography}
\end{document}